\documentclass[12pt,a4paper]{article}
\usepackage{amssymb,amsmath,theorem}
 \let\kappa=\varkappa                  
\let\emptyset=\varnothing                                  
\newcommand\norm[1]{\left\|#1\right\|}                     
\newcommand\qed{\ifhmode\unskip\nobreak\fi\quad            
   \ifmmode\square\else\hbox{$\square$}\fi}                
\newcommand\pin{\kern.0833em}                              

\theorempreskipamount=\medskipamount \theorempostskipamount=\medskipamount

 \newcommand\proof[1]{{\it Proof#1}\,}
 \newtheorem{theorem}{Theorem}
 \newtheorem{lemma}[theorem]{Lemma}
 \newtheorem{corollary}[theorem]{Corollary}
 {\theorembodyfont{\normalfont}
 \newtheorem{remark}{Remark}
 \newtheorem{example}{Example}}

\begin{document}

\begin{center}
 {\bf \large Analysis of relationships between spectral potential of transfer \\[3pt] operators,  $\boldsymbol t$-entropy, entropy and topological pressure}

\bigskip\medskip\normalsize\rm

V.\,I.\ BAKHTIN$^*$

\smallskip

{\it Belarusian State University, Belarus $($e-mail: bakhtin@tut.by\/$^*$$)$}

\bigskip

A.\,V.\ LEBEDEV

\smallskip

{\it Belarusian State University, Belarus $($e-mail: lebedev@bsu.by\/$)$}
\end{center}

\vspace{-18pt}

\renewcommand\abstractname{}
\begin{abstract} \noindent
The paper is devoted to the analysis of relationships between principal objects of the spectral theory of dynamical systems (transfer and weighted shift operators) and basic characteristics of information theory and thermodynamic formalism (entropy and topological pressure). We present explicit formulae linking these objects with $t$-entropy and spectral potential. Herewith we uncover the role of inverse rami-rate, forward entropy along with essential set and the property of non-contractibility of a dynamical system.
\end{abstract}

\bigbreak\bigskip

\quad\parbox{13.9cm} {\textbf{Keywords:} {\itshape spectral potential, transfer operator, entropy, topological pressure, $t$-entropy, essential set, rami-rate, forward entropy}

\medbreak

}

\bigskip

 \tableofcontents

\section{Introduction}

Transfer and weighted shift operators are the principal objects as in the theory of dynamical
systems so also in numerous fields of analysis. One cannot also overestimate the role of entropy
and topological pressure as in information theory so also in foundations of thermodynamic
formalism.

The paper is devoted to the analysis of interrelations between the spectral radii of the mentioned
operators and topological pressure, entropy, $t$-entropy and arising herewith dynamical and metric
invariants.

\smallskip

Let $X$ be a Hausdorff compact space, $C(X)$ be the Banach space of continuous functions on~$X$
equipped with the uniform norm, and $\alpha\!: X\to X$ be a continuous mapping. This mapping
generates  a dynamical system with discrete time which will be denoted by $(X,\alpha)$.

A linear operator $A\!: C(X)\to C(X)$ is called a \emph{transfer operator} for the dynamical system
$(X,\alpha)$ if

a) it is positive  (i.\,e., it maps nonnegative functions to nonnegative ones) and

b) it satisfies the \emph{homological identity}
\begin{equation}
\label{e-homolog}
 A\bigl((f\circ \alpha)\cdot g\bigr) = fAg, \qquad f,g\in C(X).
\end{equation}

A typical (popular) example of a transfer operator is the classical Perron--Frobenius
operator of the form
\begin{equation} \label{2,,1}
A f(x) \pin := \sum_{y\in\alpha^{-1}(x)} {a (y)}f(y),
\end{equation}

 \medskip\noindent
where $a \in C(X)$ is a certain nonnegative function. This operator is well defined when $\alpha$
is a local homeomorphism and acts onto.

If $\alpha$ is a homeomorphism then transfer operator turns out to be a classical weighted shift
(weighted composition) operator
\begin{equation} \label{2,,1a}
 Af(x) = a(x)f\bigr(\alpha^{-1}(x)\bigr),
\end{equation}
where $a\in C(X)$ and $a\ge 0$.

Transfer and weighted shift operators have numerous applications in dynamical systems theory, mathematical physics and in particular in thermodynamics, stochastic processes, information theory, investigations of zeta functions, Fredholm determinants, and operator algebras theory. They serve as an inexhaustible source of important examples and counterexamples and so also as key constructive elements of the crossed product algebras, in the theory of solvability of functional differential equations, wavelet analysis, etc. We refer the reader to the books \cite{Anton,Anton_Lebed,ABL98,ABL98', Bal00,KS97,PP90,PU10,Rue78,Rue91}, recent papers
\cite{ABL11',ABL12,ABLS03,BK19,Did07, Ex03,Kit99,Kwa12,KL13,Rou96}, and the bibliography therein.

Given a transfer operator $A$ we define a family of operators $A_\psi\!: C(X) \to C(X)$ depending
on the functional parameter $\psi \in C(X,{\mathbb R})$, where $C(X,{\mathbb R})$ is the space of continuous  real-valued
functions, by means of the formula
\begin{equation} \label{e-A-psi}
 A_\psi f := A(e^\psi f).
\end{equation}
Evidently, all the operators of this family are transfer operators as well. Let us denote
by~$\lambda(\psi )$ the logarithm of the spectral radius of $A_\psi $, that is
\begin{equation} \label{e-lamb-psi}
 \lambda(\psi ) = \lim_{n\to\infty}\frac{1}{n}\ln \norm{A_\psi^n}.
\end{equation}
The functional $\lambda(\psi )$ is called \emph{spectral potential} of $A$.

Spectral properties of weighted shift and transfer operators and especially the formulae and
methods for calculation of their spectral radii are tightly related to the ergodic and entropy
theory of dynamical systems via variational principles of thermodynamic and informational nature.

Namely, for the case when $\alpha$ is a homeomorphism and a transfer operator $A$ is a shift
operator $Af(x) =f(\alpha^{-1}(x))$ the variational principle for spectral potential $\lambda
(\psi)$ (i.\,e., variational principle for the spectral radius of weighted shift operator) has the
form
$$
 \lambda (\psi) \pin= \max_{\mu \in M_\alpha(X)} \mu[\psi] \pin=
 \max_{\mu\in E\!\pin M_\alpha(X)} \mu[\psi],
$$
where $M_\alpha(X)$ is the set of all Borel $\alpha$-invariant probability measures and $E\!\pin
M_\alpha(X)$ is the set of all ergodic measures on $X$, and
$$
 \mu[\psi] := \int_X \psi \, d\mu \pin.
$$
This variational principle was established independently by Kitover \cite{Kitover} and Lebedev
\cite{Lebedev79}. Its applications to elliptic theory of functional differential operators and
spectral theory of operator algebras associated with automorphisms are presented in \cite{Anton},
\cite{Anton_Lebed} and \cite{ABL98}. A comprehensive analysis of the corresponding variational
principles and their interrelations as with integrals so also with Lyapunov exponents for abstract
weighted shift operators associated with endomorphisms of Banach algebras is presented in
\cite{KL20}.

Spectral analysis of Perron--Frobenius operators (that is, transfer operators arising in the
situation when $\alpha$ is a local homeomorphism) naturally involves an additional dynamical object,
namely, topological pressure $P(\alpha, \psi), \ \ \psi \in C(X, {\mathbb R})$ (a detailed definition of topological pressure is
given in Subsection~\ref{s1-1}). Topological pressure appears, in particular, in the analysis of
complexity of dynamical systems $(X,\alpha)$, where $X$ is a compact metric space, and it is also a
principal component of thermodynamic formalism.

Ruelle--Walters variational principle \cite{{Rue73},Rue89,Wal75,Wal82} expresses the topological
pressure as
$$
P(\alpha, \psi)= \sup_{\mu\in M_\alpha(X)} \left(\mu[\psi] + h_\alpha(\mu)\right),
$$
where $h_\alpha(\mu)$ is Kolmogorov--Sinaj entropy.

Let $X$ be a compact metric space, $\alpha\!: X\to X$ be a local homeomorphism, and
\begin{equation*}
 Af(x) \pin= \sum_{y\in\alpha^{-1}(x)} f(y)
\end{equation*}
be the initial transfer operator in $C(X)$. Then one of the fundamental principles of thermodynamic
formalism can be written as
\begin{equation} \label{e-intr-spectr-TP}
 \lambda(\psi) = P(\alpha, \psi) \ \ \textrm{for an expanding map} \ \, \alpha.
\end{equation}
A number of results in this direction is known, cf. \cite{Bow75}, \cite{Rue78}, \cite{Wal78},
\cite{LS88}, \cite{Rue89}, \cite{LM98}, \cite{FJ01}, \cite{PU10}. However, none of these sources
considers a general case: usually it is assumed that $\alpha$ is topologically mixing, $e^\psi$ is
Holder continuous, and the space $X$ is a finite dimensional manifold or a shift space. Recently in
\cite{BK19} it is proven that $\lambda(\psi) = P(\alpha, \psi)$ for an arbitrary open expanding map
$\alpha\!: X \to X$ on a compact metric space and arbitrary  function $\psi\in C(X, {\mathbb R})$.

Topological pressure itself is a weighted version of topological entropy $h(\alpha)= P(\alpha, 0)$,
which is one of the principal ingredients of information theory. Here the corresponding variational
principle was established by Dinaburg~\cite{Din70} and Goodman~\cite{Good71}:
$$
h(\alpha) = \sup_{\mu\in M_\alpha(X)} h_\alpha(\mu).
$$

The variational principle for spectral potential of a weighted shift and general transfer operators (i.\,e., in the situation when $X$ is not necessarily a metric space and $\alpha \!: X\to X$ is an arbitrary continuous mapping) was established in \cite{B10}, \cite{AnBakhLeb11} and has the form
$$
 \lambda(\psi )=\max_{\mu\in M_\alpha(X)} \bigl(\mu[\psi]+\tau(\mu)\bigr),
$$
where $\tau(\mu)$ is a new entropy type object, called $t$-entropy, that depends not only on
dynamical system $(X,\alpha)$ but on a generating transfer operator $A$ too (see \cite{AnBakhLeb11}, \cite{BL2017}), and is a principal ingredient in entropy statistic theorem \cite{BL19}.

The latter variational principle means, in particular, that spectral potential $\lambda (\psi)$ is
the Legendre transform of $-\tau(\mu)$ and this observation also naturally leads to its
applications in thermodynamic formalism \cite{ABLS03}, \cite{ABL12}.

\smallskip

We have to emphasize that the aforementioned objects: spectral potential, topological pressure, entropy, $t$-entropy have different nature and origination and in general they can not be reduced to each other (see, in particular, Example~\ref{ex1-spec-pot-TP}). At the same time the
phenomena when they do relate to each other play the principal role and clarify the internal
structural basement of the corresponding fields of analysis (cf. \eqref{e-intr-spectr-TP}).

The goal of the paper is to uncover reasons for $\lambda(\psi )$, $h(\alpha)$, $P(\alpha, \psi)$, $\tau(\mu)$, and $h_\alpha (\mu)$ to be related to each other and describe the
arising relationships.

The article is organized as follows.

The starting Section~\ref{s-1} is devoted to the introductory overview of our principal operator, spectral, and dynamical heroes along with inevitably arising technical objects and instruments. Here we recall the notions of the spectral potential, $t$-entropy, topological entropy and pressure and the corresponding variational principles (Subsections~\ref{ss1-2},~\ref{s1-1}), introduce and examine in short the set $X_\alpha$ of essential points, i.\,e.,\ the domain where all the variational principles in question live (Subsection~\ref{ss-1-ess}), describe (in Subsection~\ref{ss-11}) the relationship between transfer operators and positive functionals,
and consider in Theorems~\ref{t-tau-omega}~and~\ref{t-lambda-omega} the `traces' on $X_\alpha$ of transfer operators, $t$-entropies and spectral potentials (later on these objects play a vital role in the analysis of problems under investigation). We also show here by examples (in particular, in Examples~\ref{ex1-compat-4}~and~\ref{ex1-spec-pot-TP}) as relationships so also drastic differences between the target objects of the research.

The main part of the article starts with Section~\ref{s-2}. Its goal is to bring to light the internal dynamical and metrical  reasons for spectral potential,  topological pressure and integrals with respect to invariant measures to be related to each other. The arising relations (estimates) are  transparently presented in Theorems~\ref{t.pressure-lambda} and~\ref{t.VP-lambda} by means of  a subsidiary object --- essential spectral potential and newly introduced dynamical invariants: rami-rate $\omega (\alpha)$ and forward entropy $\gamma (\alpha)$ that also serve as a convenient  instrument for estimating and evaluating the  topological entropy (Lemmas~\ref{l.tau+ph},~\ref{phi=0} and Corollary~\ref{c.tau+ph-h}). Moreover, we uncover here the dynamical-metric properties (Property ($*$) and Property ($**$)) in the presence of which the estimates obtained become strict equalities (Theorems~\ref{ph=0-TP} and \ref{t-lambda-phi-press-ext}). In this context an important observation is revealing (in Lemma~\ref{l-exp-trans-phi}) of a wide class of dynamical systems (generated by open non-contracting mappings) possessing the mentioned properties.

In Section~\ref{s3} we return back to the principal objects of  analysis and highlight the situations when the spectral potential $\lambda (\psi)$ is equal to  topological pressure $P(\alpha, \psi +\ln \rho)$ (Theorems~\ref{t-3-1}~and~\ref{t3-2}). This is done on the base of the results obtained in the preceding section along with a thorough analysis of the properties of the cocycle $\rho$ defining the `trace' $A_{X_\alpha}$ of the transfer operator $A$ on the set  $X_\alpha$ of essential points (Lemma~\ref{l-coc}).

Up to  Section~\ref{s-non-neg} we consider mainly transfer operators $A_\psi$ from \eqref{e-A-psi} with positive weights $e^\psi, \ \psi \in C(X,\mathbb R)$. After that in Section~\ref{s-non-neg} we
pursue the theme further and analyse the arising relationships between spectral radii, topological pressure and integrals for transfer operators with nonnegative (not necessarily positive) weights.

And  finally in  Section~\ref{s-4} we show how the results obtained can be used for explicit calculation of $t$-entropy by means of integrals and Kolmogorov--Sinaj entropy.

\section[Starters: spectral potential, topological pressure, entropy,\\
    $\boldsymbol t$-entropy etc.]{Starters: spectral potential, topological pressure,\\
    entropy, $\boldsymbol t$-entropy etc.}\label{s-1}

In this section we introduce and discuss a number of objects that will be inevitable
in the analysis of problems in question.

\subsection[$T$-entropy and variational principle for spectral potential]
  {$\boldsymbol T$-entropy and variational principle for spectral potential} \label{ss1-2}

Let us start with the spectral potential $\lambda(\psi )$ \eqref{e-lamb-psi}, i.\,e., the logarithm of the spectral radius of $A_\psi $ \eqref{e-A-psi}.

The positivity of transfer operator implies that
\begin{equation} \label{2,,2}
 \lambda(\psi ) = \lim_{n\to\infty}\frac{1}{n}\ln \norm{A_\psi^n\mathbf 1},
\end{equation}
where \textbf{1} is the unit function on $X$, and $\|\,{\cdot}\,\|$ denotes the uniform norm.

The principal object related to the spectral potential is $t$-entropy.

For an $\alpha$-invariant probability measure $\mu$ (and only such measures will be essential in our considerations) $t$-entropy $ \tau(\mu)$ is defined in the following way \cite{BL2017}:
\begin{equation} \label{e-tau-n}
 \tau(\mu):=\, \pin \inf_{n\in \mathbb N}\frac{1}{n}\inf_G \sum_{g\in G}
 \mu [g]\ln\frac{\mu[A^n g]}{\mu[g]},
\end{equation}
where the infimum $\inf_G$ is taken over the set of all continuous partitions of unity $G$ in $C(X)$ and we assume that if $\mu [g]=0$ then the corresponding summand in the right hand part is
equal to $0$ independently of the value of $\mu[A^n g]$.

Note parenthetically that if one identifies a Borel measure $\mu$ on $X$ with a linear functional
$\mu\!: C(X,\mathbb R) \to \mathbb R$ given by
\begin{equation*}
 \mu[f] := \int_X f\, d\mu,
\end{equation*}

 \smallskip\noindent
then by Riesz's theorem there exists a unique regular Borel measure on $X$ defining the same
functional. Thus, since in the foregoing definition of $t$-entropy only continuous functions (forming partitions of unity) were exploited we can assume that $t$-entropy is defined namely
for regular measures $\mu$. By $M_\alpha(X)$ we denote the set of regular $\alpha$-invariant
probability measures.

Recently an essentially different formula for $t$-entropy was obtained. Namely, it is proven in \cite[Theorem 21]{BL20} that
\begin{equation} \label{e-tau-KL1}
 \tau(\mu) \pin=\pin
 \inf_{n\in \mathbb N}\,\frac{1}{n}\!\pin\int_X \ln\frac{d(A^{*n}\mu)_a}{d\mu}\,d\mu \pin,
\end{equation}
where\/ $A^*\!: C(X)^* \to C(X)^*$ is the operator adjoint to\/ $A$, and $(A^{*n}\mu)_a$ is the
absolutely continuous component of the measure $A^{*n}\mu$ in its decomposition into absolutely
continuous and singular parts with respect to $\mu$.

The relation between $t$-entropy and spectral potential is given by the next
\begin{theorem} \label{5..6}
{\rm (variational principle for spectral potential \cite[Theorem~5.6]{AnBakhLeb11})\,} Suppose
$A\!: C(X)\to C(X)$ is a transfer operator for a continuous mapping\/ $\alpha\!: X\to X$ of a
Haus\-dorff compact space\/~$X$. Then its spectral potential\/~$\lambda(\psi )$ satisfies the
variational principle
\begin{equation} \label{5,,12}
 \lambda(\psi ) \pin= \max_{\mu\in M_\alpha(X)} \bigl(\mu[\psi]+\tau(\mu)\bigr),
 \qquad\psi \in C(X,\mathbb R).
\end{equation}
\end{theorem}

\subsection{Entropy and topological pressure} \label{s1-1}

Among the principal heroes in what follows will be the topological pressure and
topological entropy.
Therefore we recall the corresponding definitions.

These objects are defined for a dynamical system $(X,\alpha)$, where $X$ is a compact metric space
with a metric $d$.

The definitions exploit the so-called $(n,\varepsilon)$-spanning and $(n, \varepsilon)$-separated
subsets of~$X$. Let us describe them.

For every $n \in \mathbb N$ we consider the metric $d_n$ on $X$ given by
$$
 d_n (x,y) := \max \big\{ d\bigl(\alpha^i (x), \alpha^i(y)\bigr) \bigm| i=0,1,\dots,n-1\big\}.
$$

For any $\varepsilon >0$ a set $E\subset X$ is called $(n,\varepsilon)$-\emph{spanning} if it is an
$\varepsilon$-net for $X$ with respect to metric $d_n$, that is for any $x\in X$ there exists $y\in
E$ such that $d_n (x,y) <\varepsilon$.

A set $F \subset X$ is called $(n, \varepsilon)$-\emph{separated} if for each pair of points $x,y
\in F$, \ $x\neq y$, one has $d_n (x, y) >\varepsilon$.

Definition of the \emph{topological entropy} is the following:
\begin{align} \label{e.top-entr}
 h(\alpha) :={}&
 \ln \lim_{\varepsilon \to 0} \varlimsup_{n\to\infty}\inf
 \bigl\{|E|^{1/n} : \text{$E$ is $(n,\varepsilon)$-spanning}\bigr\} \\[3pt]  \label{e.top-entr1}
 ={}&\ln \lim_{\varepsilon \to 0} \varlimsup_{n\to\infty}\sup
 \bigl\{|E|^{1/n} : \text{$E$ is $(n,\varepsilon)$-separated}\bigr\}.
\end{align}

Topological pressure is a (weighted) generalization of the notion of topological entropy. Namely,
for each positive function $a\in C(X,\mathbb R)$ the \emph{topological pressure} $P(\alpha, \ln a)$ is given
by the formula
\begin{align} \label{e.top-press}
 P(\alpha, \ln a) :={} &\ln \lim_{\varepsilon \to 0} \varlimsup_{n\to\infty}\inf
 \Biggl\{\Bigg(\sum_{y\in E}\prod_{i=0}^{n-1} a\big(\alpha^i (y)\big)\!\Bigg)^{\! 1/n} :\,
 \text{$E$ is $(n,\varepsilon)$-spanning}\Biggr\}\\[3pt] \label{e.top-press1}
 ={} &\ln \lim_{\varepsilon \to 0} \varlimsup_{n\to\infty}\sup \Bigg\{\Bigg(\sum_{y\in F}
 \prod_{i=0}^{n-1} a\big(\alpha^i (y)\big)\!\Bigg)^{\! 1/n} :\,
\text{$E$ is $(n,\varepsilon)$-separated}\Biggr\}.
\end{align}
Clearly,
$$
h(\alpha) = P(\alpha,0).
$$

 \medskip

The well known Dinaburg--Goodman variational principle for topological entropy
\cite{Din70,Good71}
states that
\begin{equation} \label{e.VP-entropy}
 h(\alpha) \pin= \sup_{\mu\in M_\alpha(X)} h_\alpha(\mu)\pin,
\end{equation}

 \smallskip\noindent
and the Ruelle--Walters variational principle for topological pressure \cite{{Rue73},Rue89,Wal75,Wal82} has the form
\begin{equation} \label{e.VP-TP}
 P(\alpha, \ln a) \pin= \sup_{\mu\in M_\alpha(X)} \bigl(\mu[\ln a] + h_\alpha(\mu)\bigr),
\end{equation}

 \smallskip\noindent
where $h_\alpha(\mu)$ is the metric (Kolmogorov--Sinaj) entropy.

\subsection[Essential set $X_\alpha$]{Essential set $\boldsymbol{X_\alpha}$}
\label{ss-1-ess}

In all the mentioned variational principles for topological entropy \eqref{e.VP-entropy}, spectral potential \eqref{5,,12}, and topological pressure \eqref{e.VP-TP} the invariant measures are of essential use. In this subsection we discuss the principal set associated with their supports.

A point $x\in X$ will be called \emph{essential} (for the mapping $\alpha$), iff for every its
neighborhood $V(x)$ there exists an invariant measure $\mu\in M_\alpha(X)$ such that $\mu(V(x))
>0$. Clearly the set of all essential points is closed. We will call it the \emph{essential set}
of the mapping $\alpha$ and denote by $X_\alpha$. The points $x\in X\setminus X_\alpha$ will be
called \emph{inessential}. For each $x\in X\setminus X_\alpha$ there exists a neighborhood $V(x)$
such that $\mu(V(x)) =0$ for all $\mu\in M_\alpha(X)$. This implies that for a compact metric space
$X$ the support of each invariant measure belongs to $X_\alpha$, and for an arbitrary Hausdorff
compact space $X$ the support of each regular invariant measure belongs to $X_\alpha$.

\begin{remark} \label{r-TP-omega}
Since for each $\mu \in M_\alpha(X)$ one has $\mathop{\mathrm{supp}} \mu \subset X_\alpha$ the
variational principles \eqref{e.VP-entropy} and \eqref{e.VP-TP} imply that in the definitions of
topological entropy \eqref{e.top-entr}, \eqref{e.top-entr1} and topological pressure
\eqref{e.top-press}, \eqref{e.top-press1} it is enough to confine ourselves only to
$(n,\varepsilon)$-spanning and $(n,\varepsilon)$-separated sets in $X_\alpha$.
\end{remark}

A subset $Y\subset X$ is called \emph{$\alpha$-invariant} if $\alpha^{-1}(Y) = Y$; and it is \emph{forward $\alpha$-invariant} if $\alpha(Y) \subset Y$.

\begin{lemma} \label{l-cover}
For the essential set $X_\alpha$ one has
\begin{equation} \label{e-1-2-X-al}
 \alpha(X_\alpha)= X_\alpha.
\end{equation}
\end{lemma}

\proof. Let us check first forward $\alpha$-invariance of $X_\alpha$. Suppose on contrary that
$x\in X_\alpha$ while $\alpha (x)\notin X_\alpha$. Take a neighborhood $V(\alpha (x))$ mentioned in
the definition of inessential points, i.\,e., such that
\begin{equation} \label{e-X-a}
 \mu\big(V(\alpha (x))\big) =0 \ \ \textrm{for each} \ \, \mu\in M_\alpha(X).
\end{equation}
Then $\alpha^{-1}\big(V(\alpha (x))\big)$ is a neighborhood of $x$ and $\alpha$-invariance of
$\mu$ along with \eqref{e-X-a} imply
$$
\mu\big(\alpha^{-1}\big(V(\alpha (x))\big)\big) =0 \ \ \textrm{for each} \ \, \mu\in M_\alpha(X),
$$
that contradicts the assumption $x\in X_\alpha$.

\smallskip

Now let us prove the equality $\alpha(X_\alpha)= X_\alpha$. Suppose on contrary, that
\begin{equation} \label{,,21}
 X_\alpha \setminus \alpha(X_\alpha) \neq \emptyset.
\end{equation}
Since $\alpha(X_\alpha)$ is compact the set
$$
 W:= X\setminus \alpha(X_\alpha)
$$

 \smallskip\noindent
is open, has nonempty intersection with $X_\alpha$, and
$$
 \alpha^{-1}(W) \subset X\setminus X_\alpha.
$$
As we have mentioned the support of each regular invariant measure belongs to $X_\alpha$. Thus
$$
 \mu \bigl(\alpha^{-1}(W)\bigr) = 0 \ \ \textrm{for \ each} \ \, \mu\in M_\alpha(X).
$$
Therefore by $\alpha$-invariance of $\mu$ one has
$$
 \mu(W)= \mu\bigl(\alpha^{-1}(W)\bigr) = 0 \ \ \textrm{for \ each} \ \, \mu\in M_\alpha(X),
$$
that contradicts \eqref{,,21}. \qed

\begin{remark} \label{r11-X-al-not-invar}
In spite of the equality $X_\alpha = \alpha (X_\alpha)$, the set of essential points $X_\alpha$ in
general  is not $\alpha$-invariant, i.\,e., it may occur that $\alpha^{-1}(X_\alpha) \neq
X_\alpha$. The next example demonstrates such a phenomenon.
\end{remark}

\begin{example} \label{e11-X-al-not-invar}
Let $X=[0,1]\subset \mathbb R$ and
$$
 \alpha(x)=
   \begin{cases}
     2x, & x\in [0,1/2],\\[2pt]
     1, & x\in[1/2,1].
   \end{cases}
$$
Routine check shows that there are only two ergodic measures for this $(X,\alpha)$, namely, the
Dirac measures $\delta_0$ and $\delta_1$. Thus $X_\alpha = \{0,1\}$, while $\alpha^{-1}(X_\alpha) =
\{0\}\cup [1/2,1] \neq X_\alpha$.
\end{example}

The next result gives statistic criteria describing the essential set.

Recall that for every $x\in X$ the \emph{empirical measure} $\delta_{x,n}$ is defined by the formula
\begin{equation*}
 \delta_{x,n} =\frac{1}{n}\big( \delta_x +\delta_{\alpha(x)} +\,\dotsm\,+\delta_{\alpha^{n-1}(x)} \big).
\end{equation*}

\begin{theorem} \label{essential}
The following three conditions are equivalent:

\smallskip

$a)$ a point\/ $x\in X$ is essential;

\smallskip

$b)$ for any neighborhood\/ $V(x)$ of\/ $x$ there exists\/ $y\in X$ such that
\begin{equation} \label{,,liminf}
 \liminf_{n\to \infty} \delta_{y,n}\big(V(x)\big) >0;
\end{equation}

$c)$ for any neighborhood\/ $V(x)$ of\/ $x$ there exists\/ $y\in X$ such that
\begin{equation} \label{,,limsup}
 \limsup_{n\to \infty} \delta_{y,n}\big(V(x)\big) >0.
\end{equation}
\end{theorem}

\proof. For any point $x$ and any its neighborhood $V(x)$ there exists a neighborhood $U(x)$ such that its closure belongs to $V(x)$. By  Urysohn's lemma there exists a function $f\in C(X)$ that is equal to $1$ on $U(x)$, equal to $0$ outside $V(x)$, and takes values in $[0,1]$ on $V(x)\setminus U(x)$.

$a)\ \Rightarrow\ b)$. If a point $x$ is essential then there is a measure $\mu\in M_\alpha(X)$ such that $\mu\big(U(x)\big) >0$. Therefore $\mu[f]>0$. By the choice of $f$ and the Ergodic theorem \cite[\S 1.6]{Wal82} there exists a measurable function $\bar f\!:X\to [0,1]$ satisfying the conditions
\begin{equation*}
 \bar f =\bar f\circ\alpha, \quad \mu\big[\bar f\pin\big] =\mu[f], \quad \text{and}\quad
 \delta_{y,n}[f]\to \bar f(y) \ \ \text{for $\mu$-almost all $y$}.
\end{equation*}
This implies
\begin{equation} \label{,,liminf_greater}
 \liminf_{n\to\infty} \delta_{y,n}\big(V(x)\big) \ge \liminf_{n\to\infty} \delta_{y,n}[f]
 =\bar f(y) \ \ \text{a.e.}
\end{equation}

 \medskip\noindent
Since $\mu[\bar f] =\mu[f] >0$ it follows that  the function $\bar f(y)$ is positive on a set of positive measure $\mu$. And if $\bar f(y)>0$ then \eqref{,,liminf_greater} implies \eqref{,,liminf}.

$b)\ \Rightarrow\ c)$. Obvious.

$c)\ \Rightarrow\ a)$. Let $V(x), U(x)$  be the neighbourhoods and the function $f$ mentioned above. By virtue of \eqref{,,limsup} with $U(x)$ substituted for $V(x)$  there exists a point $y\in X$, a number $\varepsilon>0$, and an infinite subset $\mathbb N_\varepsilon \subset \mathbb N$ such that
\begin{equation} \label{,,delta_greater}
 \delta_{y,n}[f] \ge \delta_{y,n}\big(U(x)\big) >\varepsilon \ \ \text{for all $n\in \mathbb N_\varepsilon$}.
\end{equation}
By the Banach--Alaoglu theorem the set $\{\pin \delta_{y,n}\mid n\in\mathbb N_\varepsilon\pin\}$ possesses a limit point $\mu$ in the dual space $C(X)^*$ equipped with the $^*$-weak topology. This $\mu$ is a linear functional on $C(X)$, which by Riesz's theorem is identified with a regular probability measure on $X$. In a standard manner it can by verified that this measure is $\alpha$-invariant, i.e., $\mu\in M_\alpha(X)$. Passing to the limit in \eqref{,,delta_greater} we obtain $\mu\big(V(x)\big) \ge \mu[f] \ge\varepsilon$. By the arbitrariness of $V(x)$ it follows that $x$ is essential. \qed

\medskip

Recall that a point $x \in X$ is \emph{non-wandering} if for every open neighborhood $V$ of $x$ we
have $V \cap \alpha^n(V)\neq \emptyset$ for some $n \in \mathbb N$. The set of non-wandering points
is denoted by $\Omega(\alpha)$. It is a closed forward $\alpha$-invariant set \cite[Theorem~5.6.]{Wal82}. We have $\text{supp}\, \mu \subset \Omega (\alpha)$ for every $\mu \in M_\alpha(X)$.

Clearly wandering points are inessential and therefore $X_\alpha\subset \Omega(\alpha)$. In reality
the latter inclusion may be strict (i.\,e., not all the non-wandering points are essential,
$X_\alpha \neq \Omega(\alpha)$). This phenomenon is demonstrated by the following example.

\begin{example} \label{ex-X-a-Omega-a}
Let us consider the standard one-third Cantor set
\begin{equation*}
 C \,:=\, \biggl\{ x=\sum_{i=1}^\infty \frac{x_i}{3^i}\biggm| x_i\in \{0,2\}\biggr\}.
\end{equation*}
It is naturally identified with the set of all sequences of the form $x =(x_1,x_2,x_3,\dots)\in
\{0,2\}^{\mathbb N}$. Consider the mapping $\alpha\!: C\to C$, \ $\alpha(x) :=3x \pmod 1$. This
mapping is the left shift on the set of sequences. Let $X\subset C$ be the set consisting of all
the sequences $(x_1,x_2,\dots)\in \{0,2\}^{\mathbb N}$ such that each segment $(x_i,\dots,
x_{i+2^n-1})$ of length $2^n$ contains not more than $n$ copies of digit $2$ (for each
$i,n\in\mathbb N$).

By routine check one sees that $\alpha(X) =X$, the set $X$ is closed and all its points are
non-wandering. Note also that for the mapping $\alpha \!: X\to X$ there is a unique essential point $x^* =(0,0,0,\dots)$, i.\,e., $X_\alpha =\{x^*\}$. To prove this it is enough to verify that the
dynamical system $(X,\alpha)$ possesses a unique ergodic measure $\delta_{x^*}$.

Let us consider an arbitrary ergodic measure $\mu$ for $(X,\alpha)$. Take a finite sequence
$y=(y_1,\dots,y_m)\in \{0,2\}^m$. It defines a (may be empty) cylinder
\begin{equation*}
 Z_m(y) \,=\, \big\{ x=(x_1,x_2,\dots)\in X\bigm| x_1=y_1,\,\dots,\,x_m=y_m\big\}.
\end{equation*}
By the ergodic theorem for $\mu$-almost all points $x\in X$ the relative number of points of the
trajectory $x$, $\alpha(x)$, \dots, $\alpha^{2^n-1}(x)$ got into $Z_m(y)$ converges to
$\mu(Z_m(y))$. On the other hand, if $y=(y_1,\dots,y_m)$ contains at least one `two', then by
definition of $X$ this relative number does not exceed $m(n+1)/2^n$. Therefore, $\mu(Z_m(y)) =0$ for each cylinder $Z_m(y)$ that do not contain $x^*$ and thus the measure $\mu$ is supported at $x^*$.
\end{example}

In the analysis of asymptotic properties of trajectories of dynamical systems it is often enough to
consider the restriction of the initial mapping $\alpha$ onto the set of non-wandering points
$\Omega(\alpha)$. In contrast, in our paper we consistently adhere to the point of view that in the
analysis of invariant measures and characteristics associated with them (such as entropy and
topological pressure) it is enough to consider the restriction of the initial mapping $\alpha$ onto
the essential set $X_\alpha$.

\begin{remark} \label{r-omega-inf}
For the set $\Omega(\alpha)$ of non-wandering points one has $\alpha(\Omega(\alpha)) \subset \Omega(\alpha)$. Thus one can consider the non-wandering set  $\Omega_2(\alpha):= \Omega(\alpha|_{\Omega(\alpha)})\subset \Omega(\alpha)$. It can occur that $\Omega_2(\alpha)\neq \Omega(\alpha)$ (see, for example,  \cite[\S~5.3]{Wal82}). If we put $\Omega_1(\alpha):=\Omega(\alpha)$ and define inductively $\Omega_n(\alpha):=\Omega(\alpha|_{\Omega_{n-1}(\alpha)}), \ n=2,3, \dots$ then $\Omega_1(\alpha)\supset \Omega_2(\alpha)\supset \cdots$ is a decreasing set of closed subsets
of $X$. So one can put
$$
\Omega_\infty(\alpha):=\bigcap_{n=1}^\infty\Omega_n(\alpha) .
$$
For the same reasons as for  $\Omega(\alpha)$ we have  that $\mu(\Omega_\infty(\alpha))=1$ for each $\mu \in M_\alpha(X)$.

One can also consider other subsets of  $\Omega(\alpha)$ possessing the mentioned property. For example, we can take the set $R(\alpha)\subset \Omega(\alpha)$ of the so-called \emph{recurrent points} (see, \cite[\S~6.4]{Wal82}) that also possesses the property: $\mu(R(\alpha))=1$ for each $\mu \in M_\alpha(X)$.

Clearly, the set $X_\alpha$ is the minimal one possessing this property.
\end{remark}

\subsection{Transfer operators and positive functionals. Compatibility} \label{ss-11}

In this subsection we present a procedure of `taking traces' of transfer operators that will play an important technical role in the further analysis.

We start with  a  more explicit description of transfer operators linking them with
special families of positive functionals.

Let, as above, $X$ be a Hausdorff compact space, $\alpha\!: X\to X$ be a continuous mapping, and $A\!: C(X)\to C(X)$ be a certain transfer operator.

For every point $x\in X$ define the functional $\phi_x$ according to the formula
\begin{equation} \label{7,,1}
 \phi_x [f] := \bigl[{A}f\bigr](x), \qquad f\in C(X).
\end{equation}
In other words,
\begin{equation} \label{7,,1a}
 \phi_x := A^* \delta_x,
\end{equation}

 \medskip\noindent
where $A^*\!: C(X)^* \to C(X)^*$ is the adjoint to $A$ operator and $\delta_x$ is the Dirac
functional
$$
 \delta_x [f] = f(x), \qquad f\in C(X).
$$
Evidently, $\phi_x$ is a positive functional.

There are two possibilities for $x$.

\smallskip

1) \,$[A{\mathbf 1}] (x)=0$. This means that $\phi_x [{\mathbf 1}]=0$ which implies $\phi_x =0$ due
to the positivity of $\phi_x$.

\smallskip

2) \,$[A{\mathbf 1}](x)\neq 0$. In this case $\phi_x \neq 0$ and $\phi_x$ defines a regular
measure $\nu_x =A^*\delta_x$ on $X$.

\smallskip

The homological identity \eqref{e-homolog} implies also that for any $f \in C(X)$ we have
$$
 \bigl[A(f\circ\alpha)\bigr](x) =\bigl[A\bigl((f\circ\alpha)\cdot{\mathbf 1}\bigr)\bigr](x) =
 f(x)\cdot [A{\mathbf 1}](x),
$$
and therefore
$$
 \frac{1}{[A{\mathbf 1}](x)}\,\phi_x (f \circ \alpha)= f (x),
$$
which means that
\begin{equation} \label{7,,2}
 \mathop{\mathrm{supp}} \nu_x \subset \alpha^{-1} (x).
\end{equation}

\medskip

Clearly, the mapping $x \mapsto \phi_x$ is $^*$-weakly continuous on $X$.


The family $\{\phi_x\}$ presented above in fact gives a complete description of transfer operators in $C(X)$ since one can easily verify that every $^*$-weakly continuous mapping $x \mapsto \phi_x$, where $\phi_x$ are positive functionals satisfying \eqref{7,,2} (here $\phi_x$ may be $0$ as well), defines a certain transfer operator $A\!: C(X)\to C(X)$ acting according to formula~\eqref{7,,1}.


%

\begin{remark} \label{r-loc-hom-Perron}
The foregoing description of transfer operators implies, in particular, that in the situation when
$\alpha \!: X\to X$ is a local homeomorphism each transfer operator $A$ acts as classical
Perron--Frobenius operator \eqref{2,,1} on $\alpha (X)$ and $Af\vert_{X\setminus \alpha (X)}
\equiv 0$ for any  $f\in C(X)$; and in the situation when $\alpha \!: X\to X$ is a homeomorphism each transfer operator is a weighted shift operator \eqref{2,,1a}.
\end{remark}

For any closed subset $Y\subset X$ we can naturally define the `trace' $\phi_{Y,x}$ of the
functional $\phi_x$ from~\eqref{7,,1} on $C(Y)$. Here is its definition. We can identify each
function $f\in C(Y)$ with the function $\tilde{f}$ on $X$ of the form
\begin{equation} \label{e-tilde-f}
 \tilde{f}(x) :=
  \begin{cases}
     f(x), & x\in Y, \\[2pt]
     0, & x\notin Y.
  \end{cases}
\end{equation}
 Since each positive functional on $C(X)$ is defined by a unique regular Borel measure we can
 uniquely extend its values onto the functions of the form \eqref{e-tilde-f} and in this way for
 the aforementioned functionals $\phi_x$ we set
\begin{equation} \label{e-transf-Y0}
 \phi_{Y,x}(f) := \phi_x\big(\tilde{f}\pin\big).
\end{equation}

\smallskip

The next notion is inevitable for taking `traces' of transfer operators.

\smallskip

Let  $Y\subset X$ be a  closed forward $\alpha$-invariant set.  A transfer operator $A\!: C(X) \to C(X)$  will be called \emph{$Y$-compatible} (or \emph{compatible with $Y$}) iff the family of functionals $\phi_{Y,x}$ on $C(Y)$
defined by \eqref{e-transf-Y0} is $^*$-weakly continuous on $Y$.

To clarify the notion introduced we present examples as of $Y$-compatible so also not
$Y$-compatible operators.

\begin{example} \label{ex1-compat-1}
Let $Y$ be any closed $\alpha$-invariant set, then any transfer operator $A$ is $Y$-compatible. This follows from definition \eqref{e-transf-Y0} of $\phi_{Y,x}$ and property \eqref{7,,2} of
$\phi_{x}$.
\end{example}



\begin{example} \label{ex1-compat-3}
Even in the case when $\alpha$ is a homeomorphism and $Y$ is a  closed forward $\alpha$-invariant set there can exist transfer operators that are not $Y$-compatible.

Let $X=[0,1]$ and  $\alpha (x)=x^2$. If $A: C(X)\to C(X)$ is a transfer operator then by the foregoing description we have $Af(x) = \rho (\sqrt{x})f(\sqrt{x})$, where $\rho \in C(X), \ \rho \geq 0$. That is $\phi_x = \rho (\sqrt{x}) \delta_{\sqrt{x}}$\pin. Set $Y=[0,x_0]$, where $0<x_0 <1$. This $Y$ is a closed forward $\alpha$-invariant set. And we have that
$$
 \phi_{Y,x}=
  \begin{cases}
    \phi_x, & \sqrt{x}\leq x_0, \\[2pt]
    0, & \sqrt{x}> x_0.
  \end{cases}
$$

 \medskip\noindent
That is $A$ is $Y$-compatible iff $\rho(x_0)=0$.
\end{example}

\begin{example} \label{ex1-compat-4}
Even in the situation when $Y$ is closed and $\alpha(Y)=Y$ there can exist transfer operators that are not $Y$-compatible.

Let
\begin{equation} \label{e-Bakht-X}
 X= [0,1]\times [0,1] = \Delta_1 \sqcup \Delta_2,
\end{equation}
where
\begin{gather*}
 \Delta_1:= \biggl\{(x_1,x_2)\in X \biggm| x_2\le \frac{2-x_1}{2}\biggr\}, \\[9pt]
\Delta_2:= \biggl\{(x_1,x_2)\in X \biggm| x_2 > \frac{2-x_1}{2}\biggr\}.
\end{gather*}

\smallskip

Set
\begin{equation}\label{e-Bakh}
 \alpha (x_1,x_2) \,=\,
\begin{cases}
  \Big(x_1, \sqrt{2x_2/(2-x_1)}\pin\Big), & (x_1,x_2)\in  \Delta_1 , \\[6pt]
  (x_1,1), & (x_1,x_2)\in  \Delta_2 .
\end{cases}
\end{equation}
And let
\begin{equation}\label{e-Bakht-A}
 [Af](x_1,x_2) = f\big(x_1, x_2^2 (2-x_1)/2\big), \quad (x_1,x_2)\in X.
\end{equation}

 \medskip\noindent
Thus we have
$$
\phi_{(x_1,x_2)} =\delta_{(x_1,\, x_2^2 (2-x_1)/2)}, \quad (x_1,x_2)\in X.
$$

 \medskip\noindent
Take $Y := [0,1]\times \{0,1\}$. We have $\alpha(Y)=Y$ and
$$
\phi_{Y,(x_1,\, x_2)}=
\begin{cases}
  \delta_{(x_1,\,x_2)} , & (x_1,x_2)\in  [0,1]\times \{0\}, \\[3pt]
  \delta_{(0,1)}, & (x_1,x_2) =(0,1), \\[3pt]
  0, & (x_1,x_2)\in (0,1] \times \{1\} .
  \end{cases}
$$
Therefore $A$ is not $Y$-compatible.

On the other hand one can verify in a routine way that for a given transfer operator $A$ associated with the mapping \eqref{e-Bakh} and the mentioned set  $Y$ the operator $A$ is $Y$-compatible  if $\phi_{(0,1)}=0$, where $\phi_{(x_1,\,x_2)}$, \ $(x_1,x_2)=x \in X$, is the family of functionals \eqref{7,,1a}.
\end{example}

It is worth mentioning here that though in Example~\ref{e11-X-al-not-invar} $X_\alpha$ is not $\alpha$-invariant it is a discrete set and therefore any transfer operator $A$ is $X_\alpha$-compatible.  However, in general a transfer operator $A$ is not necessarily $X_\alpha$-compatible. This possibility is demonstrated by the next example.

\begin{example} \label{ex1-compat-10}
Let us consider the objects  mentioned in Example~\ref{ex1-compat-4}, i.\,e., $X$ \eqref{e-Bakht-X}, $\alpha$ \eqref{e-Bakh} and $A$ \eqref{e-Bakht-A}. One can check in a routine way that here we have $X_\alpha =\Omega (\alpha) = [0,1]\times \{0,1\}$. And we have already verified in Example~\ref{ex1-compat-4} that $A$  is not compatible with this set.
\end{example}

The next two examples present popular situations when $X_\alpha$ is such that any   transfer operator $A$ is $X_\alpha$-compatible.

\begin{example} \label{e.invert-A-compat}
Let $\alpha : X\to X$ be a homeomorphism. Since $\alpha(X_\alpha)= X_\alpha$ it follows that in this case $X_\alpha$ is $\alpha$-invariant and therefore any transfer operator~$A$  is $X_\alpha$-compatible.
\end{example}

\begin{example} \label{e.expand-A-compat}
Let $\alpha : X\to X$ be a local homeomorphism and $Y\subset X$ be a closed subset such that
$\alpha (Y)= Y$. If $\alpha$ is a local homeomorphism on $Y$ then observation in Remark~\ref{r-loc-hom-Perron} implies that any transfer operator $A$ is $Y$-compatible.

In particular, if $(X,d)$ is a compact metric space, $\alpha\!: X\to X$ is an \emph{expanding} map (i.\,e., the map for which there exist $r >0$ and $\Lambda>1$ such that inequality $d(x,y) \leq r$ implies $d(\alpha(x),\alpha(y))\geq \Lambda\pin d(x,y)$), and if, additionally, $\alpha$ is an open map then we have $\Omega(\alpha) = \overline{\textrm{Per} (\alpha)}$, where $\textrm{Per}(\alpha)$ is the set of periodic points
\cite[Proposition~3.3.6]{PU10}. Thus in this case we have $X_\alpha=\Omega (\alpha)
=\overline{\textrm{Per}(\alpha)}$. Moreover, in this situation $\alpha|_{\overline{\mathrm{Per}
(\alpha)}}$ is an open map as well \cite[Lemma~3.3.10]{PU10}. Therefore observation in Remark
\ref{r-loc-hom-Perron} implies that any transfer operator $A$ is $X_\alpha$-compatible.
\end{example}

Unfortunately, it can occur that the restriction of a local homeomorphism $\alpha$ onto $Y$ is not a local homeomorphism. In this case it can happen that a transfer operator $A$ is not $Y$-compatible. Here is an example of such situation.

\begin{example} \label{e.loc-hom-not-A-comp}
Let $X= [0,1]\times \Delta$, where $\Delta = \{0\}\cup \{\pin 1/2^n\mid n=0,1,2,\dots\}$, and topology on $X$ is induced from ${\mathbb R}^2$. Define $\alpha : X\to X$ by the formulae
\begin{align*}
  \alpha(t,0) = (t,0),& \qquad t\in [0,1]; \\[3pt]
  \alpha(t,{1}/{2^n}) = \left(t,{1}/{2^{n-1}}\right),& \qquad  t\in [0,1], \quad n=1,2,\dots; \\[3pt]
  \alpha (t,1) = \big(\sqrt{t},1\big),& \qquad   t\in [0,1] .
\end{align*}
Clearly $\alpha$ is a local homeomorphism.

Take $Y=[1/2,1]\times \Delta$. We have that $Y$ is a closed set, $\alpha(Y)=Y$ while $\alpha :Y\to Y$ is not a local homeomorphism at the point $(1/2,1)$.

According to observation in Remark~\ref{r-loc-hom-Perron} any transfer operator  $A :C(X)\to C(X)$ is of the form \eqref{2,,1}. Therefore $A$ is $Y$-compatible iff $a(1/2,1) =0$.
\end{example}

Given a dynamical system $(X,\alpha)$, a transfer operator $A : C(X)\to C(X)$, and a set $Y$ such that $A$ is $Y$-compatible  one can define a transfer operator $A_Y \!: C(Y) \to C(Y)$ for the
dynamical system $(Y,\alpha)$, that can be naturally considered as the `trace' of $A$ on $C(Y)$.
Namely, we set
\begin{equation} \label{e-transf-Y}
 \bigl[{A_Y}f\bigr](x):=\phi_{Y,x}[f], \qquad f\in C(Y), \quad x\in Y
\end{equation}

 \medskip\noindent
(cf. \eqref{7,,1} and \eqref{e-transf-Y0}). The argument exploited for $A$ and $\phi_x$ proves also
that $A_Y$ is a transfer operator for $\alpha\!: Y\to Y$.

\begin{remark} \label{r1-trnsf-not-compat}
Note that in the situation when $Y$ is a closed and forward $\alpha$-invariant set the mapping
$A_Y$ given by \eqref{e-transf-Y} is defined on $C(Y)$ and is `nearly' a transfer operator: it is
positive by positivity of $\phi_{Y,x}$ and satisfies by construction the homological identity, but
if $A$ is not $Y$-compatible (i.\,e., the family $\phi_{Y,x}$ is not $^*$-weakly continuous on $Y$)
then $A_Y$ does not preserve $C(Y)$. The foregoing Examples~\ref{ex1-compat-3}, \ref{ex1-compat-4} and \ref{e.loc-hom-not-A-comp} can also be considered as illustrations of  such situations.
\end{remark}



Once a dynamical system $(X,\alpha)$ and a transfer operator $A$ are fixed then each pair $(Y,A_Y)$, consisting of a   set $Y\subset X$ such that $A$ is $Y$-compatible and the above described transfer operator $A_Y$~\eqref{e-transf-Y}, defines a $t$-entropy $\tau_Y(\mu)$ on the set $M_\alpha (Y)$ of $\alpha$-invariant probability measures for the dynamical system
$(Y,\alpha)$. This $\tau_Y(\mu)$ is given by formula \eqref{e-tau-n} with $A_Y$ substituted for
$A$. Note also that each measure $\mu\in M_\alpha (Y)$ can be considered as $\mu \in M_\alpha(X)$
by setting $\mu (X\setminus Y):=0$ and in this way we assume that $M_\alpha (Y) \subset
M_\alpha(X)$.

\begin{theorem} \label{t-tau-omega}
Let\/ $A$ be a transfer operator for a dynamical system\/ $(X,\alpha)$ and\/ $Y\subset X$ be a subset such that $A$ is $Y$-compatible. Then we have
\begin{equation} \label{e1-tY}
 \tau (\mu) = \tau_{Y} (\mu), \qquad \mu \in M_\alpha(Y),
\end{equation}
and, if\/ $A$ is\/ $X_\alpha$-compatible, then
\begin{equation} \label{e1-t-Om}
 \tau (\mu) = \tau_{X_\alpha} (\mu), \qquad \mu \in M_\alpha(X).
\end{equation}
\end{theorem}

\proof. It is enough to prove \eqref{e1-tY}, since as we have noted for each $\mu\in M_\alpha(X)$
one has $\mathop{\mathrm{supp}} \mu \subset X_\alpha$ and so $M_\alpha(X) =M_\alpha(X_\alpha)$.


Recall that we are identifying measures $\mu \in M_\alpha(Y)$ with the measures $\mu\in
M_\alpha(X)$ such that $\mathop{\mathrm{supp}}\mu\subset Y$. By \eqref{e-tau-KL1} we have
\begin{equation} \label{e-tau-KL}
 \tau(\mu) \pin=\pin
 \inf_{n\in \mathbb N}\,\frac{1}{n}\!\pin\int_X \ln\frac{d(A^{*n}\mu)_a}{d\mu}\,d\mu \pin,
\end{equation}

 \medskip\noindent
where\/ $A^*\!: C(X)^* \to C(X)^*$ is the operator adjoint to\/ $A$, and $(A^{*n}\mu)_a$ is the
absolutely continuous component of the measure $A^{*n}\mu$ in its decomposition into absolutely
continuous and singular parts with respect to $\mu$.

By the construction of $A_{Y}$ (cf. \eqref{e-tilde-f}, \eqref{e-transf-Y0}, \eqref{e-transf-Y}) and
the already mentioned inclusion $\text{supp}\, \mu \subset Y$ we have that
$$
(A^{*n}\mu)_a = (A_{Y}^{*n}\mu)_a\pin .
$$
This along with \eqref{e-tau-KL} implies
$$
\tau(\mu) \pin=\pin
 \inf_{n\in \mathbb N}\,\frac{1}{n}\!\pin\int_X \ln\frac{d(A_{Y}^{*n}\mu)_a}{d\mu}\,d\mu \pin =
 \tau_{Y} (\mu), \qquad \mu\in M_\alpha(Y),
$$
that proves \eqref{e1-tY}. \qed

\medskip

Let $Y\subset X$ and $A$ be any $Y$-compatible transfer operator. For the transfer operator $A_Y \!: C(Y) \to C(Y)$ we denote by
$\lambda_Y (\psi)$ the spectral potential of operator $A_Y$, i.\,e., given by formula \eqref{2,,2}
with $A_Y$ substituted for $A$ and exploiting restriction of $\psi$ onto $Y$. Inclusion
$M_\alpha(Y)\subset M_\alpha(X)$ along with Theorems~\ref{5..6} and \ref{t-tau-omega} imply the
next

\begin{theorem} \label{t-lambda-omega}
Let\/ $A\!: C(X)\to C(X)$ be a transfer operator for a dynamical system\/ $(X,\alpha)$ and\/
$Y\subset X$ be a subset such that $A$ is\/ $Y$-compatible. Then
\begin{equation} \label{e1-lam-Y}
 \lambda(\psi )\geq\lambda_{Y} (\psi) = \max_{\mu\in M_\alpha(Y)} \bigl(\mu[\psi]+
 \tau_{Y}(\mu)\bigr),\qquad \psi \in C(X,\mathbb R);
\end{equation}
and, if\/ $A$ is\/ $X_\alpha$-compatible, then
\begin{equation} \label{e1-lam-Omega}
 \lambda(\psi ) = \lambda_{X_\alpha} (\psi) = \max_{\mu\in M_\alpha(X)} \bigl(\mu[\psi]+
 \tau_{X_\alpha}(\mu)\bigr),\qquad \psi \in C(X,\mathbb R).
\end{equation}
\end{theorem}

\begin{remark} \label{r1-lambda-lambda-Y}
Note also that the inequality $\lambda(\psi )\geq\lambda_{Y} (\psi) $ follows directly from the
relationship between $A$ and $A_Y$ (cf. \eqref{7,,1} and \eqref{e-transf-Y}) and \eqref{2,,2}.
\end{remark}

We finish the section with an example demonstrating that
spectral potential $\lambda(\psi)$ and $t$-entropy $\tau(\mu)$ cannot be reduced to topological pressure $P(\alpha,\psi)$ and entropy $h_\alpha(\mu)$, respectively.

\begin{example} \label{ex1-spec-pot-TP}
Let us consider once more the objects from Example~\ref{ex1-compat-4}, i.\,e., $X$ \eqref{e-Bakht-X}, $\alpha$ \eqref{e-Bakh} and $A$ \eqref{e-Bakht-A}. Clearly here we have $\Omega (\alpha) = [0,1]\times \{0,1\}$ and a measure $\mu$ is  $\alpha$-invariant iff $\textrm{supp}\, \mu \subset \Omega (\alpha)$, i.\,e., $X_\alpha = \Omega (\alpha)$.

For Kolmogorov--Sinaj entropy one has
\begin{equation} \label{e-KSE-Bakh}
 h_\alpha(\mu) =0, \qquad \mu\in M_\alpha(X) .
\end{equation}
And therefore  according to variational principle \eqref{e.VP-TP} for a real-valued continuous function $\psi$ the topological pressure $ P(\alpha, \psi)$ is equal to
\begin{equation}\label{e-TP-Bakh}
 P(\alpha, \psi)=\sup_{\mu\in M_\alpha(X)} \mu[\psi] = \max_{x\in X_\alpha} \, \psi (x).
\end{equation}
On the other hand formula \eqref{e-tau-KL} implies that for $t$-entropy we have
\begin{equation}\label{e-T-entr-Bakh}
 \tau(\mu)=
 \begin{cases}
   0, & \textrm{supp}\, \mu \subset \Delta , \\[2pt]
   -\infty, & \textrm{supp}\, \mu \cap (X_\alpha\setminus \Delta) \neq \emptyset ,
 \end{cases}
\end{equation}
where
$$
 \Delta = [0,1]\times \{0\}\cup \{(0,1)\},
$$

 \medskip\noindent
and therefore according to variational principle \eqref{5,,12} for the spectral potential $\lambda (\psi)$ one obtains
\begin{equation} \label{e-lam-Bakh}
  \lambda (\psi) = \max_{x\in \Delta } \, \psi(x) .
\end{equation}
\end{example}

\medskip

Comparing formulae \eqref{e-KSE-Bakh} and \eqref{e-T-entr-Bakh} one concludes that Kolmogorov--Sinaj entropy and $t$-entropy are essentially different objects. And comparing formulae \eqref{e-TP-Bakh} and \eqref{e-lam-Bakh} one arrives at the same conclusion for topological pressure and spectral exponent.

\medskip

Moreover, we have to emphasize that while Kolmogorov--Sinaj entropy $h_\alpha (\mu)$ is always non-negative the $t$-entropy $\tau (\mu)$ can take negative and even infinite negative values. And, what is important, these infinite values of $\tau (\mu)$ are quite natural:  according to variational principle \eqref{5,,12} they indicate the measures that do not play any role in the spectral potential calculation.

Now after we have introduced the key heroes of our study and established their principal difference
we are going to uncover the analytic reasons for them to be related to each other and describe
these relationships.

\section[Essential spectral potential, rami-rate, forward entropy]{Essential spectral potential, rami-rate, forward\\ entropy} \label{s-2}

We start with introduction of a number of characteristics of dynamical systems that will be of use
in estimation of entropy, spectral potential and topological pressure.

Given a dynamical system $(X,\alpha)$,
we put
\begin{equation} \label{e2-0}
 \tilde{\alpha}^{-n} (x):= \alpha^{-n}(x) \cap X_\alpha, \qquad x \in X_\alpha.
\end{equation}

Recall that $\alpha(X_\alpha) =X_\alpha$ by Lemma \ref{l-cover}. Henceforth we will always assume that $\alpha$ is a \emph{finite-sheeted cover on} $X_\alpha$, i.\,e., satisfies the condition
\begin{equation} \label{e.preim}
 \sup_{x\in X_\alpha} |\tilde{\alpha}^{-1}(x)| < \infty.
\end{equation}

In what follows we need in two more notions.

The number
\begin{equation} \label{e.tau}
 \omega (\alpha) := \ln \varlimsup_{n\to\infty} \sup_{x\in X_\alpha} |\tilde{\alpha}^{-n}(x)|^{1/n}
\end{equation}

 \medskip\noindent
will be called the \emph{inverse rami-rate}. It evaluates the ramification speed of $\alpha$ preimages.

In the case when $X$ is a compact metric space we put
\begin{equation} \label{e.ph}
 \gamma (\alpha) := \ln \lim_{\varepsilon \to 0} \varlimsup_{n\to\infty}\inf
 \bigl\{|\alpha^n(E)|^{1/n} : \text{$E$ is $(n,\varepsilon)$-spanning for $X_\alpha$}\bigr\}.
\end{equation}
Comparing formulae \eqref{e.top-entr} and \eqref{e.ph} we naturally call $\gamma(\alpha)$ \emph{forward entropy}.

\smallskip

Note that since $|\alpha^{n}(E)| \leq |E|$ we have (by definitions \eqref{e.top-entr}, \eqref{e.ph}
along with Remark~\ref{r-TP-omega}) that
\begin{equation} \label{e.ph-entr-neq}
 \gamma (\alpha) \leq h(\alpha).
\end{equation}

 \medskip\noindent
There are examples when $\gamma (\alpha) < h(\alpha)$ (see, in particular, Lemma \ref{phi=0}) so in general forward entropy $\gamma (\alpha)$ and topological entropy $h(\alpha)$ are different characteristics
of~$\alpha$.

Relationship between inverse rami-rate $\omega (\alpha)$ and topological entropy $h(\alpha)$ is based on the following

\begin{lemma} \label{l.tau-h}
Let\/ $X$ be a compact metric space and\/ $\alpha \!: X\to X$ be a local homeomorphism. Then there
exists\/ $\varepsilon > 0$ such that for each\/ $n\in \mathbb N$ and\/ $x\in X$ the set\/
$\alpha^{-n} (x)$ is\/ $(n,\varepsilon)$-se\-par\-ated.
\end{lemma}

\proof. Since $\alpha$ is a local homeomorphism it follows that $|\alpha^{-1}(x)|$ is a continuous
(locally constant) function. The set $ \Delta := \big\{x\in X : |\alpha^{-1}(x)| \le 1\big\}$ is
clopen in $X$. Thus the set $X\setminus \Delta$ is compact.

Let us put $D \!: X\setminus \Delta \to (0,\infty)$ to be
$$
 D(x):= \min\bigl\{d(u,v)\bigm| u,v\in \alpha^{-1}(x),\ u\ne v\bigr\}.
$$
Local homeomorphness of $\alpha$ implies that $D(x)$ is a continuous function. Therefore
$$
\varepsilon := \frac{1}{2}\min \bigl\{ D(x)\bigm| x\in X\setminus \Delta\bigr\} >0.
$$
Routine verification shows that this $\varepsilon$ fits the statement of lemma.
\qed

\medskip

This lemma along with \eqref{e.top-entr1} and Remark~\ref{r-TP-omega} implies that in the situation
when $\alpha$ is a local homeomorphism on $X_\alpha$ one has
\begin{equation} \label{e.neq-tau-h}
 \omega(\alpha) \leq h(\alpha).
\end{equation}

\begin{remark} \label{r-h-t}
Since, for example, for any invertible $\alpha$ we have $\omega(\alpha)=0$, while (by choosing
suitable invertible $\alpha$) \,$h(\alpha)$ could be any nonnegative number (see, for example, \cite[\S 7.3]{Wal82}), we conclude that $\omega(\alpha)$ and $h(\alpha)$ are different characteristics of $\alpha$.
\end{remark}

\begin{lemma} \label{l.tau+ph}
The characteristics\/ $h(\alpha)$, $\gamma (\alpha)$, and\/ $\omega(\alpha)$ satisfy the inequality
\begin{equation} \label{e.tau+ph}
 h(\alpha) \le \gamma (\alpha) +\omega(\alpha).
\end{equation}
\end{lemma}

\proof. Recalling Remark~\ref{r-TP-omega} we obtain
\begin{align*}
 h(\alpha) &=
 \ln \lim_{\varepsilon \to 0} \varlimsup_{n\to\infty}\inf \bigl\{|E|^{1/n} : \text{$E$ is
 $(n,\varepsilon)$-spanning for $X_\alpha$}\bigr\}\\[3pt]
 &\le
 \ln \lim_{\varepsilon \to 0} \varlimsup_{n\to\infty}\Bigl[\inf\bigl\{|\alpha^n(E)|^{1/n} :
 \text{$E$ is $(n,\varepsilon)$-spanning for $X_\alpha$}\bigr\}\times
 \sup_{x\in X_\alpha}|\tilde{\alpha}^{-n}(x)|^{1/n}\Bigr]\\[3pt]
  &\leq
 \ln \Bigl[\lim_{\varepsilon \to 0} \varlimsup_{n\to\infty}\inf\bigl\{|\alpha^n(E)|^{1/n} :
 \text{$E$ is $(n,\varepsilon)$-span. for $X_\alpha$}\bigr\}\times \varlimsup_{n\to\infty}
 \sup_{x\in X_\alpha}|\tilde{\alpha}^{-n}(x)|^{1/n}\Bigr]\\[3pt]
 &= \gamma(\alpha) + \omega(\alpha). \qed
\end{align*}

This lemma along with observation \eqref{e.ph-entr-neq} implies

\begin{corollary} \label{c.tau+ph-h}
If inequality\/ \eqref{e.neq-tau-h} holds then

\smallskip

$(i)$ \,if\/ $\omega (\alpha) =0$ then\/ $h(\alpha) =\gamma (\alpha)$,

\smallskip

$(ii)$ \,if\/ $\gamma (\alpha) =0$ then\/ $h(\alpha) =\omega (\alpha)$.
\end{corollary}

\begin{remark}
Inequality \eqref{e.tau+ph} may be strict and equalities $h(\alpha) =\gamma(\alpha)$ and
$h(\alpha) =\omega(\alpha)$ may take place not only in the case when the second
summand (i.\,e., $\omega(\alpha)$ or $\gamma(\alpha)$, respectively) is zero (see
Example~\ref{ex.tau+ph}).
\end{remark}

In what follows we will make use of the next auxiliary spectral potential type object which will
help us to present the results in a transparent way.

Let $(X,\alpha)$ be a dynamical system with $\alpha$ being a finite-sheeted cover on $X_\alpha$. For each nonnegative function $a\in C(X)$ we put
\begin{equation} \label{e.sp-exp}
 \ell (\alpha, a) := \ln \lim_{n\to \infty} \sup_{x\in X_\alpha} \Biggl(\,
 \sum_{y\in\tilde{\alpha}^{-n}(x)} \prod_{i=0}^{n-1} a\bigl(\alpha^i (y)\bigr)\!\Biggr)^{\! 1/n},
\end{equation}
where we set $\ln 0 = -\infty$. The number $\ell (\alpha, a)$ will be called
\emph{essential spectral potential}.

\begin{remark} \label{r2-sp-pot-al}
Note that\/ $\ell(\alpha, a)$ is the logarithm of the `spectral radius' of Perron--Fro\-be\-nius
operator\/ $A$ associated with\/ $(X_\alpha, \alpha)$, i.\,e., with the dynamical system on the essential set $X_\alpha$ (and that is why we use the term essential spectral potential). We put here the `spectral radius' in quotation marks since in general $($when\/ $\alpha$ is not a local homeomorphism\/$)$ formula\/ \eqref{2,,1} does not define an
operator in\/ $C(X)$.
\end{remark}

The next result links topological pressure with essential spectral potential via forward entropy $\gamma(\alpha)$.

\begin{theorem} \label{t.pressure-lambda}
Let\/ $X$ be a compact metric space, $\alpha \!: X\to X$ be a local homeomorphism on\/ $X_\alpha$, and\/ $a\in C(X)$ be a positive function. Then
$$
P(\alpha, \ln a) - \gamma (\alpha)\leq \ell (\alpha, a) \leq P(\alpha, \ln a).
$$
\end{theorem}

\proof. The right-hand inequality follows from \eqref{e.sp-exp}, \eqref{e.top-press1} along with
Lemma~\ref{l.tau-h} and Remark~\ref{r-TP-omega}.

To prove the left-hand inequality note that for each finite subset $E\subset X_\alpha$ one has
\begin{align} \label{e.P-L}
 \sum_{y\in E}\prod_{i=0}^{n-1} a\bigl(\alpha^i(y)\bigr) &\leq
 \sum_{x\in \alpha^n(E)}\sum_{y\in \tilde{\alpha}^{-n}(x)}\prod_{i=0}^{n-1} a\bigl(\alpha^i(y)\bigr)\\
 \label{e.P-L1}
 & \leq |\alpha^n(E)|\times \sup_{x\in X_\alpha}\sum_{y\in \tilde{\alpha}^{-n}(x)}\prod_{i=0}^{n-1}
 a\bigl(\alpha^i(y)\bigr).
\end{align}
Denoting for brevity
\begin{equation} \label{F-sup}
 \Phi_n:= \sup_{x\in X_\alpha}\sum_{y\in \tilde{\alpha}^{-n}(x)}
 \prod_{i=0}^{n-1} a\bigl(\alpha^i(y)\bigr),
\end{equation}

 \smallskip\noindent
we can rewrite formula \eqref{e.sp-exp} in the form
\begin{equation} \label{e.lambda0}
 \ell(\alpha,a) =\lim_{n\to\infty}\frac{1}{n}\ln\Phi_n.
\end{equation}
Now observation \eqref{e.P-L}, \eqref{e.P-L1} along with formula \eqref{e.top-press} for $P(\alpha,
\ln a)$, Remark~\ref{r-TP-omega} and definition \eqref{e.ph} of $\gamma (\alpha)$ implies
\begin{align*}
 P(\alpha, \ln a) &\leq \lim_{\varepsilon \to 0}
 \varlimsup_{n\to\infty} \frac{1}{n}\ln\Bigl[\inf\bigl\{|\alpha^n(E)| :
 \text{$E$ is $(n,\varepsilon)$-spanning for $X_\alpha$}\bigr\}\times \Phi_n\Bigr]\\[3pt]
 &= \lim_{\varepsilon \to 0}  \varlimsup_{n\to\infty} \frac{1}{n}\ln\inf\bigl\{|\alpha^n(E)| :
 \text{$E$ is $(n,\varepsilon)$-spanning for $X_\alpha$}\bigr\} +\lim_{n\to\infty}\frac{1}{n}
 \ln \Phi_n\\[3pt]
 &= \gamma (\alpha) + \ell(\alpha, a). \qed
\end{align*}

\begin{corollary} \label{TP-phi=0}
Under conditions of Theorem \ref{t.pressure-lambda} if\/ $\gamma (\alpha) =0$ then\/ $P(\alpha,
\ln a) = \ell(\alpha, a)$.
\end{corollary}

\begin{remark} \label{r-P-neq-l-phi-neq-0}
If $\gamma (\alpha) >0$ then we can have
$$
P(\alpha, \ln a) - \gamma (\alpha)= \ell(\alpha, a) < P(\alpha, \ln a).
$$
Indeed, let $\alpha\!: X\to X$ be a homeomorphism. Then $\omega(\alpha)=0$ and $h(\alpha)
=\gamma(\alpha)$. By a suitable choice of $X$ and $\alpha$ one can assume that $h(\alpha)
=\gamma(\alpha)$ is an arbitrary given nonnegative number. For this $X$ and $\alpha$ take also
$a=1$. Then $\ell(\alpha, 1) =0$ and
$$
 P(\alpha, \ln 1) =P(\alpha, 0) = h(\alpha) = \gamma(\alpha) > 0.
$$
\end{remark}

Theorem~\ref{t.pressure-lambda} shows importance of forward entropy $\gamma (\alpha)$. This
characteristics can be easily calculated in the presence of the next

 \medskip\noindent
\textbf{Property} $\boldsymbol{(*)}$ \ For each pair $(n,\varepsilon)$, \ $n\in\mathbb N$, \,$\varepsilon
>0$, there exists a finite set $F(n,\varepsilon)\subset X_\alpha$ such that
the set $\tilde{\alpha}^{-n}\bigl(F(n,\varepsilon)\bigr)$ is an $(n,\varepsilon)$-spanning for
$X_\alpha$ and $\lim_{n\to\infty} |F(n,\varepsilon)|^{1/n} =1$.

\medskip

This property looks as being rather sophisticated. A particular (more convenient) variant is the next

\medskip

\textbf{Property} $\boldsymbol{(**)}$ \ For each $\varepsilon >0$ there exists a finite set
$F(\varepsilon)\subset X_\alpha$ such that for each $n\in \mathbb N$ the set $\tilde{\alpha}^{-n
}(F(\varepsilon))$ is an $(n,\varepsilon)$-spanning for $X_\alpha$.

\medskip

Clearly Property $(**)$ implies Property $(*)$ since one can simply take $F(n,\varepsilon) :=
F(\varepsilon)$ for all $n\in \mathbb N$.

\begin{lemma}  \label{phi=0}
If\/ $\alpha$ possesses property\/ $(*)$ then\/ $\gamma(\alpha) =0$ {\rm(}and hence\/ $h(\alpha) =
\omega(\alpha)${\rm)}.
\end{lemma}

\proof. By definition \eqref{e.ph} of $\gamma (\alpha)$
$$
 \gamma(\alpha) \pin\leq\pin \ln\lim_{\varepsilon \to 0} \varlimsup_{n\to\infty}
 \bigl|\alpha^n\bigl(\tilde{\alpha}^{-n}(F (n,\varepsilon))\bigr)\bigr|^{1/n} \pin=\pin
 \ln\lim_{\varepsilon \to 0} \varlimsup_{n\to\infty} \bigl|F(n,\varepsilon)\bigr|^{1/n} \pin=\pin 0.  \qed
$$

\smallskip

As a consequence of Lemma~\ref{phi=0} and Theorem~\ref{t.pressure-lambda} we also obtain

\begin{theorem} \label{ph=0-TP}
If\/ $\alpha$ possesses property\/ $(*)$ then\/ $P(\alpha, \ln a) = \ell(\alpha, a)$.
\end{theorem}

Lemma~\ref{l-exp-trans-phi} below presents a wide class of dynamical systems possessing Property
$(**)$ (and therefore Property $(*)$).

Recall that a mapping $\alpha\!: X\to X$ on a metric space $(X,d)$ is called \emph{non-contracting}
if there exists $r>0$ such that inequality $d(x,y)\le r$ implies $d(\alpha(x),\alpha(y)) \ge
d(x,y)$.

In the proof of Lemma \ref{l-exp-trans-phi} we will use the following technical observation.

\begin{lemma} \label{shadowing}
Let\/ $X$ be a compact metric space and\/ $\alpha\!:X\to X$ be a non-contracting local homeomprphism. Then there exists an\/ $\varepsilon>0$ such that inequality\/
$d(x,\alpha(y)) <\varepsilon$ implies existence of a point\/ $z\in\alpha^{-1}(x)$ such that\/
$d(z,y) \le d(x,\alpha(y))$.
\end{lemma}

\proof. Let $r$ be the number from definition of non-contractiveness of $\alpha$. By the openness
of $\alpha$ for each point $x\in X$ there exists $\varepsilon(x)>0$ such that
\begin{equation} \label{e2-compr}
 \alpha\big(B(x,r/2)\big) \supset B\big(\alpha(x),2\varepsilon(x)\big).
\end{equation}
Now for each $x\in X$ let us take a (small) neighborhood $U(x)$ such that
$$
U(x)\subset B(x,r/2) \quad \text{and} \quad
\alpha\bigl(U(x)\bigr)\subset B\bigl(\alpha(x),\varepsilon(x)\bigr).
$$
By the choice of $U(x)$ along with \eqref{e2-compr} for each point $y\in U(x)$ we have
\begin{equation} \label{open_mapping_1}
 \alpha\big(B(y,r)\big) \supset \alpha\big(B(x,r/2)\big) \supset
 B\big(\alpha(x),2\varepsilon(x)\big) \supset B\big(\alpha(y),\varepsilon(x)\big).
\end{equation}

For the family of the mentioned neighborhoods $U(x)$ there exists a finite subcover $U(x_1)$,
\dots, $U(x_n)$ of the space $X$. Set $\varepsilon :=\min\{\pin \varepsilon(x_i)\mid
i=1,\dots,n\pin\}$. Now \eqref{open_mapping_1} implies
\begin{equation} \label{open_mapping_2}
 \alpha\big(B(y,r)\big) \supset B\big(\alpha(y),\varepsilon\big), \qquad y\in X.
\end{equation}

Finally note that if $d(x,\alpha(y)) <\varepsilon$ then by \eqref{open_mapping_2} there exists
$z\in \alpha^{-1}(x) \cap B(y,r)$. And since $\alpha$ is non-contracting one has $d(z,y) \le
d(x,\alpha(y))$. \qed

\begin{lemma} \label{l-exp-trans-phi}
If the mapping\/ $\alpha\!: X\to X$ is a non-contracting local homeomorphism on\/ $X_\alpha$ then it possesses property\/ $(**)$.
\end{lemma}

\proof. It suffice to take $\varepsilon>0$ for which the statement of Lemma~\ref{shadowing} holds
and as $F(\varepsilon)$ one can take any $\varepsilon$-net in $X_\alpha$. Indeed, for each $y\in
X_\alpha$ there exists a point $x_n\in F(\varepsilon)$ such that $d(x_n,\alpha^n(y))< \varepsilon$
and by means of Lemma~\ref{shadowing} one can construct a sequence of points $x_{n-1}$,
$x_{n-2}$\pin, \dots, $x_1$, $x_0$ in $X_\alpha$ such that for all $i=1,\dots,n$ the following
conditions hold
\begin{equation*}
 x_{i-1}\in \alpha^{-1}(x_i), \qquad d\big(x_{i-1},\alpha^{i-1}(y)\big) \le \pin
 d\big(x_{i},\alpha^{i}(y)\big) <\varepsilon.
\end{equation*}

 \medskip\noindent
These relations show that the set $\alpha^{-n}(F(\varepsilon))$ forms an $(n,\varepsilon)$-spanning
in $X_\alpha$. \qed

\medskip

Summarising Lemma~\ref{phi=0}, Theorem~\ref{ph=0-TP} and Lemma~\ref{l-exp-trans-phi} we obtain

\begin{theorem} \label{t-lambda-phi-press-ext}
Suppose\/ $\alpha \!: X\to X$ is a non-contracting local homeomorphism on\/ $X_\alpha$. Then

\smallskip

$(i)$ \ $\gamma(\alpha) =0$ and thus\/ $h(\alpha) =\omega(\alpha)$$;$

\smallskip

$(ii)$ \ $P(\alpha,\ln a) = \ell(\alpha, a)$.
\end{theorem}

\begin{remark}\label{r-exp-diff}
Properties $(i)$ and $(ii)$ for expanding diffeomorphisms of compact smooth manifolds where stated
without proofs in \cite{LM98}.
\end{remark}

The next example shows that inequality in \eqref{e.tau+ph} may be strict.

\begin{example} \label{ex.tau+ph}
Let $X= S^1 \sqcup Y$, where $S^1$ is the unit circle $S^1 = \{\pin z\in\mathbb C: |z|=1\pin\}$ and
$Y$ is a certain compact metric space. Set $\alpha := \alpha_1 \oplus \alpha_2$, where
$\alpha_1(z) =z^N$, \ $z\in S^1$, and $\alpha_2$ is a homeomorphism of $Y$. By
Theorem~\ref{t-lambda-phi-press-ext} $(i)$ \ $\gamma (\alpha_1) =0$ and
$h(\alpha_1)=\omega(\alpha_1) =\ln N$ (where the latter equality follows directly from
\eqref{e.tau}). On the other hand, since $\alpha_2$ is a homeomorphism we have that $\omega
(\alpha_2) =0$, and therefore $\gamma (\alpha_2) = h(\alpha_2)$. By a suitable choice of $Y$ and
$\alpha_2$ one can assume that $\gamma (\alpha_2) = h(\alpha_2)$ is an arbitrary given nonnegative
number. Note also that $\omega (\alpha) = \omega(\alpha_1)$ and $\gamma (\alpha) = \gamma
(\alpha_2)$. Now we have
$$
 h(\alpha) = \max\{h(\alpha_1),h(\alpha_2)\}= \max \{\omega(\alpha_1),\gamma(\alpha_2)\} =
 \max \{\omega(\alpha),\gamma(\alpha)\}.
$$
In particular, when $\gamma(\alpha_2) > 0$ we have
$$
 h(\alpha) < \omega(\alpha) + \gamma (\alpha).
$$
We note in addition that one can have here $h(\alpha) = \omega(\alpha)$ or $h(\alpha) = \gamma
(\alpha)$ along with $\gamma(\alpha)\neq 0$ and $\omega(\alpha)\neq 0$.
\end{example}

We finish this section with an estimate of the  essential spectral potential by means of integrals
and inverse rami-rate $\omega (\alpha)$. Note that the next theorem is valid for an \emph{arbitrary} compact space $X$ and \emph{arbitrary} $\alpha$, i.\,e., not necessarily satisfying condition of
Theorem~\ref{t-lambda-phi-press-ext}.

\begin{theorem} \label{t.VP-lambda}
Let\/ $X$ be a compact space, $\alpha$ be a finite-sheeted cover on\/ $X_\alpha$, and\/ $a\in C(X)$ be a nonnegative function. Then
$$
 \max_{\mu\in M_\alpha(X)}\int_X \ln a\, d\mu \,=\max_{\mu\in E\!\pin M_\alpha(X)}\int_X \ln a\,
 d\mu \,\leq\, \ell(\alpha,a) \,\leq\,
 \max_{\mu\in E\!\pin M_\alpha(X)}\int_X \ln a\, d\mu\pin +\omega(\alpha).
$$
\end{theorem}

\proof. Keeping in mind that for each $\mu \in M_\alpha(X)$ we have $\textrm{supp}\pin\mu \subset
X_\alpha$ and applying $\limsup$ Variational Principle \cite[Theorem~3.5]{KL20} one obtains
\begin{equation} \label{e.inv-VP}
 \lim_{n\to \infty} \sup_{x\in X_\alpha}\frac{1}{n} \ln\Biggl(\pin\prod_{i=0}^{n-1} a\bigl(\alpha^i
 (x)\bigr)\!\Biggr) \pin= \max_{\mu\in E\!\pin M_\alpha(X)}\int_X \ln a\, d\mu \,= \max_{\mu\in
 M_\alpha(X)}\int_X \ln a\, d\mu\pin.
\end{equation}
Therefore,
\begin{align*}
 \ell (\alpha, a) &\pin=\pin \ln\lim_{n\to \infty}\sup_{x\in X_\alpha}\Biggl(\,\sum_{y\in
 \tilde{\alpha}^{-n}(x)}\prod_{i=0}^{n-1} a\bigl(\alpha^i (y)\bigr)\!\Biggr)^{\!1/n} =\pin
 \lim_{n\to \infty} \sup_{x\in X_\alpha} {\frac{1}{n}} \ln\Biggl(\,\sum_{y\in\tilde{\alpha}^{-n}(x)}
 \prod_{i=0}^{n-1}  a\bigl(\alpha^i (y)\bigr)\!\Biggr)\\[3pt]
 &\geq\pin \lim_{n\to \infty} \sup_{y\in X_\alpha}\frac{1}{n} \ln\Biggl(\pin\prod_{i=0}^{n-1}
 a\bigl(\alpha^i(y)\bigr)\!\Biggr) \pin= \max_{\mu\in E\!\pin M_\alpha(X)}\int_X \ln a\, d\mu\pin,
\end{align*}
which proves the left-hand inequality in question.

On the other hand one has
\begin{align*}
 \ell (\alpha, a) &\pin=\pin \ln \lim_{n\to \infty} \sup_{x\in X_\alpha}\Biggl(\, \sum_{y\in
 \tilde{\alpha}^{-n}(x)}\prod_{i=0}^{n-1} a\bigl(\alpha^i (y)\bigr)\!\Biggr)^{\! 1/n}\\[3pt]
 &\pin\leq\pin \ln \varlimsup_{n\to \infty} \left\{\sup_{x\in X_\alpha} |\tilde{\alpha}^{-n}(x)|^{1/n} \sup_{y\in X_\alpha}\Biggl(\pin\prod_{i=0}^{n-1} a\bigl(\alpha^i (y)\bigr)\!\Biggr)^{\! 1/n}\right\}  \\[3pt]
 & \pin=\pin \ln \varlimsup_{n\to \infty} \sup_{x\in X_\alpha} |\tilde{\alpha}^{-n}(x)|^{1/n}
 \pin+\pin \lim_{n\to \infty} \sup_{y\in X_\alpha}{\frac{1}{n}}\ln\Biggl(\pin\prod_{i=0}^{n-1}
 a\bigl(\alpha^i (y)\bigr)\!\Biggr)\\[3pt]
 &\pin=\pin \omega(\alpha) \pin+ \max_{\mu\in E\!\pin M_\alpha(X)}\int_X \ln a\, d\mu\pin.
 \end{align*}
Here in the final  equality we again exploited \eqref{e.inv-VP}.  \qed

\medskip

As an immediate corollary we have

\begin{theorem} \label{t2-invert}
Let\/ $X$ be a compact space, $\alpha$ be a homeomorphism on\/ $X_\alpha$, and\/ $a\in C(X)$ be a
nonnegative function. Then
$$
 \ell(\alpha, a) \pin= \max_{\mu\in E\!\pin M_\alpha(X)}\int_X \ln a\, d\mu\pin.
$$
\end{theorem}

\section{Spectral potential vs topological pressure} \label{s3}

Now we return back to description of interrelation between spectral potential $\lambda (\psi)$ and
topological pressure $P(\alpha,\psi)$. Throughout this section we assume that $A$ is a transfer operator for a dynamical system $(X,\alpha)$ and $A_\psi$ and $\lambda (\psi)$ are defined by \eqref{e-A-psi} and \eqref{e-lamb-psi} respectively.

From now on we adopt the following convention. Once we use a transfer operator $A$ and the
essential set $X_\alpha$ we assume that $A$ is an $X_\alpha$-compatible and $\alpha$ is a finite-sheeted cover on $X_\alpha$, i.\,e., satisfies condition \eqref{e.preim}.

For this situation  the description of transfer operator
$A_{X_\alpha}$ is given in Subsection~\ref{ss-11} (cf. \eqref{e-tilde-f}, \eqref{e-transf-Y0},
\eqref{e-transf-Y}). It implies that
\begin{equation} \label{e-3.1}
 \bigl[A_{X_\alpha}f\bigr](x) \pin= \sum_{y\in{\tilde{\alpha}}^{-1}(x)}\rho(y)f(y), \qquad f\in
 C(X_\alpha), \ \ x\in X_\alpha,
\end{equation}
where $\rho$ is a certain nonnegative function on $X_\alpha$. This function $\rho$ is usually called a \emph{cocycle} associated with the transfer operator $A_{X_\alpha}$.

In fact the cocycle $\rho$ has rather specific properties and henceforth we proceed to describe some of them.

Throughout all the discussion of these properties (up to Corollary~\ref{l-coc2}) in order not to overload the notation we will simply write $X$ instead of $X_\alpha$, $A$ instead of $A_{X_\alpha}$, and $\alpha^{-1}$ instead of $\tilde\alpha^{-1}$. In this notation our setting looks as follows: we consider a continuous finite-sheeted cover $\alpha\!:X\to X$, i.\,e., satisfying the condition
\begin{equation}\label{e-rho-infty}
   \sup_{x\in X} |\alpha^{-1}(x)| < \infty ;
\end{equation}
(cf. \eqref{e.preim}); and a transfer operator $A\!: C(X)\to C(X)$ of the form
\begin{equation} \label{e-rho-A}
 [Af](x) \pin= \sum_{y\in\alpha^{-1}(x)}\rho(y)f(y), \qquad f\in
 C(X), \ \ x\in X,
\end{equation}
where $\rho$ is a certain nonnegative function (cocycle) on $X$ (cf. \eqref{e-3.1}).

A point $x\in X$ will be called

--- a \emph{local injectivity point} (LIP) if there exists a neighborhood $U(x)$ such that the mapping $\alpha\!: U(x) \to X$ is injective;

--- a \emph{local openness point} (LOP) if for any neighborhood $U(x)$ its image $\alpha (U(x))$ contains some neighborhood of $\alpha (x)$;

--- a \emph{local homeomorphism point} (LHP) if $x$ and $\alpha(x)$ have $\alpha$-homeomorphic neighborhoods.

\begin{lemma} \label{l-coc}
Under condition\/ \eqref{e-rho-infty}

\smallskip

a$)$ if\/ $\rho(x_0)=0$ then\/ $\rho$ is continuous at the point\/ $x_0$$;$

\smallskip

b$)$ if\/ $\rho(x_0)\neq 0$ then\/ $\rho$ is continuous at\/ $x_0$ iff\/ $x_0$ is a LIP$;$

\smallskip

c$)$ if\/ $\rho(x_0) \neq 0$ then\/ $x_0$ is a LOP$;$

\smallskip

d$)$ if\/ $\rho (x_0)\neq 0$ then\/ $x_0$ is a LIP iff it is a LHP.
\end{lemma}

\proof. Using \eqref{e-rho-infty}, choose a neighborhood $O(x_0)$ such that
\begin{equation} \label{,,61}
 \alpha^{-1}\big(\alpha (x_0)\big)\cap O(x_0) = \{x_0\}.
\end{equation}
Choose a nonnegative function $f\in C(X)$ such that $f(x_0) =1$ and $f(x) =0$ outside~$O(x_0)$. Then by \eqref{e-rho-A} we have
\begin{equation} \label{,,62}
 [Af](\alpha(x_0)) =\rho(x_0),
\end{equation}
and
\begin{equation} \label{,,63}
 [Af](\alpha(x)) \ge \rho(x)f(x) \quad \text{for all}\ \ x\in X.
\end{equation}

\medskip

\textit{a}) If $\rho(x_0) =0$ then $[Af](\alpha(x_0)) =0$ by \eqref{,,62}. Along with continuity of $[Af](\alpha(x))$ and \eqref{,,63} this implies
\begin{equation*}
 0\le \rho(x) \le \frac{[Af](\alpha(x))}{f(x)} \longrightarrow 0 \quad\text{as}\ \ x\to x_0,
\end{equation*}
which means the continuity of $\rho(x)$ at $x_0$.

\textit{b}) Let $x_0$ be a LIP and $U(x_0)$ be a neighborhood where $\alpha$ is injective. Then
\begin{equation*}
  \alpha^{-1}\big(\alpha (x)\big)\cap U(x_0) = \{x\} \quad \textrm{for any}\ \ x\in U(x_0)
\end{equation*}
and hence by \eqref{e-rho-A}
\begin{equation*}
 [Af](\alpha(x)) =\rho(x) f(x) \quad \textrm{for any}\ \ x\in U(x_0).
\end{equation*}
Since $[Af](\alpha(x))$ is continuous and $f(x_0) =1$, this implies the continuity of $\rho(x)$ at $x_0$.

On the other hand, assume that $\rho(x)$ is continuous at $x_0$ but $x_0$ is not a LIP. Then there are pairs of disjoint points $x,\,x'$ both arbitrarily close to $x_0$ and such that $\alpha(x) =\alpha(x')$. For these pairs we have
\begin{equation*}
 \limsup_{x\to x_0}\pin [Af](\alpha(x)) \,\ge\, \limsup_{x\to x_0}
 \bigl(\rho(x)f(x) +\rho(x')f(x')\bigr) \,\ge\,  2\rho(x_0).
\end{equation*}
Since $[Af](\alpha(x))$ is continuous and $\rho(x_0)\ne 0$, this contradicts \eqref{,,62}.

\textit{c}) Take an arbitrary neighborhood $U(x_0)$ and a nonnegative function $g\in C(X)$ such that $g(x_0) =1$ and $g(x) =0$ outside $U(x_0)$. Then
\begin{equation*}
 [Ag](\alpha(x_0)) \ge \rho(x_0)g(x_0) >0.
\end{equation*}
Take a neighborhood $V$ of $\alpha(x_0)$ of the form
\begin{equation*}
 V:= \big\{ y\in X \bigm| [Ag](y) > 0\big\}.
\end{equation*}
Now, the choice of $g$ and definition \eqref{e-rho-A} imply $V\subset \alpha(U(x_0))$.


\textit{d}) Suppose $x_0$ is a LIP. Then by \textit{b}) the function $\rho$ is continuous at $x_0$ and hence it is positive in a certain neighborhood of $x_0$. Take a neighborhood $U(x_0)$ such that $\alpha$ is injective and $\rho$ is positive on it. By \textit{c}) all points of $U(x_0)$ are LOPs. Consequently, $\alpha$ is open on $U(x_0)$ and maps it homeomorphically onto $\alpha(U(x_0))$. \qed

\begin{remark} \label{r-LIP=LHP}
Lemma just proven demonstrates  that for a cocycle $\rho$ the property  of being continuous is
valid only in rather specific situations. Fortunately, as assertion \textit{c}) tells, at the points where $\rho$ does not vanish the mapping $\alpha$ behaves not `too pathologically'.

In general LIP and LHP are different notions. For example, the point $(1/2,1)$ of  the set $Y$ in Example~\ref{e.loc-hom-not-A-comp} is a LIP (for $\alpha\!: Y\to Y$) but is not a LHP. Lemma \ref{l-coc} shows, in addition, that a transfer operator $A$ from \eqref{e-rho-A} is a `clever machine' --- it \emph{distinguishes} LIP and LHP: when $x_0$ is a LIP but not a LHP it puts $\rho (x_0)=0$.
\end{remark}

As an immediate consequence  of the forgoing Lemma we also obtain

\begin{corollary}
\label{l-coc2}
If\/ $\rho(x_0)\neq 0$ then\/ $\rho$ is continuous at a point\/ $x_0$ iff\/ $x_0$ is a LHP.
\end{corollary}

And, in particular, for the objects mentioned in \eqref{e-3.1} one  has the following

\begin{corollary} \label{c-coc3}
If\/ $\alpha\!: X \to X$ is a local homeomorphism on\/ $X_\alpha$ then the cocycle\/ $\rho$ defined in\/ \eqref{e-3.1} is a continuous function.
\end{corollary}

Recall that we assume throughout the rest of the article that $A$ is an $X_\alpha$-compatible transfer operator and $\alpha$ is a finite-sheeted cover on $X_\alpha$.

\smallskip

The next observation links spectral potential $\lambda (\psi)$ defined in \eqref{e-lamb-psi} and
essential spectral potential $\ell(\alpha,a)$ defined in \eqref{e.sp-exp}.

\begin{theorem} \label{t-3-1a}
Let the cocycle\/ $\rho$ of\/ $A$ be continuous on\/ $X_\alpha$ $($in particular, this is true when
$\alpha \!: X_\alpha\to X_\alpha$ is a local homeomorphism\/$)$. Then
$$
 \lambda(\psi) = \ell\bigl(\alpha,\rho e^\psi\bigr),
$$

 \medskip\noindent
$($recall that\/ $\ell(\alpha,a)$ in \eqref{e.sp-exp} exploits only the values of\/ $a$
on $X_\alpha$$)$.
\end{theorem}

\proof. By \eqref{e-A-psi} we have
$$
 A_{X_\alpha, \psi}f := A_{X_\alpha}\bigl(e^\psi f\bigr), \qquad f\in C(X_\alpha).
$$
Let us denote
\begin{equation} \label{e3-3}
 a:= \rho\pin e^\psi.
\end{equation}

 \medskip\noindent
Since $A_{X_\alpha, \psi}^n$ is a positive operator we have that $\bigl\| A_{X_\alpha,
\psi}^n\bigr\| = \bigl\| A_{X_\alpha, \psi}^n{\bf 1}\bigr\|$ and routine computation shows that
\begin{equation} \label{e3-4}
 \bigl\| A_{X_\alpha, \psi}^n\bigr\| \pin=\pin \max_{x\in X_\alpha}\Biggl(\, \sum_{y\in
 {\tilde{\alpha}}^{-n}(x)} \prod_{i=0}^{n-1} a\bigl(\alpha^i (y)\bigr)\!\pin\Biggr),
\end{equation}
and therefore
\begin{equation*}
 \lambda_{X_\alpha}(\psi) \pin=\pin \ln\lim_{n\to\infty}
 \bigl\| A_{X_\alpha,\psi}^n\bigr\|^{1/n} =\pin \ln\lim_{n\to\infty} \max_{x\in X_\alpha}
 \Biggl(\, \sum_{y\in {\tilde{\alpha}}^{-n}(x)} \prod_{i=0}^{n-1} a\bigl(\alpha^i(y)\bigr)
 \!\pin\Biggr)^{\! 1/n} \!=\pin \ell (\alpha, a),
\end{equation*}
where the final equality follows from definition \eqref{e.sp-exp}.

This observation along with Theorem~\ref{t-lambda-omega} implies
$$
 \lambda(\psi) = \lambda_{X_\alpha}(\psi)=\ell (\alpha, \rho e^\psi). \qed
$$

The foregoing theorem along with results of Section~\ref{s-2} gives us a possibility to relate
spectral potential and topological pressure. This is the theme of the next

\begin{theorem} \label{t-3-1}
Let\/ $X$ be a compact metric space and\/ $\alpha\!: X\to X$ be a local homeomorphism on\/ $X_\alpha$. If\/ $\alpha$ possesses property\/ $(*)$, and the cocycle\/ $\rho$ is strictly positive on\/ $X_\alpha$ then for each strictly positive continuous extension of\/ $\rho$ onto\/ $X$ we have
$$
 \lambda(\psi) = P(\alpha,\psi +\ln\rho).
$$
\end{theorem}

\proof. Define $a$ by \eqref{e3-3}. By Theorem~\ref{t-3-1a} along with
Theorem~\ref{t.pressure-lambda} we have
$$
 \lambda(\psi) = \ell(\alpha, a) = P(\alpha, \ln a) = P(\alpha,\psi +\ln\rho).  \qed
$$

Exploiting in the foregoing proof Theorem~\ref{t-lambda-phi-press-ext} in place of
Theorem~\ref{t.pressure-lambda} one gets

\begin{theorem} \label{t3-2}
Let\/ $X$ be a compact metric space, $\alpha\!: X\to X$ be a non-contracting local homeomorphism on\/ $X_\alpha$, and the cocycle\/ $\rho$ be strictly positive on\/ $X_\alpha$. Then for each strictly positive continuous extension\/ of\/ $\rho$ onto\/ $X$ we have
$$
 \lambda(\psi) = P(\alpha,\psi+\ln\rho).
$$
\end{theorem}

\section[Spectral radii of transfer operators with nonnegative weights,\\
          topological pressure and integrals]{Spectral radii of transfer operators with
          nonnegative weights, topological pressure and integrals} \label{s-non-neg}

In the preceding sections we analysed transfer operators $A_\psi =A(e^\psi\,\cdot\,)$, \ $\psi\in
C(X)$. Here the weight (i.\,e., the function $e^\psi$) is always positive. In this section we
extend the results obtained above onto transfer operators with non-negative (not necessarily
positive) weights.

Let $A$ be a fixed transfer operator for $(X,\alpha)$. We define the family of operators $Ag \!:
C(X) \to C(X)$, where $g\in C(X)$, as
\begin{equation} \label{e4a-1}
 Ag:= A(g\,\cdot\,).
\end{equation}
Clearly, if $g\geq 0$, then $Ag$ is a transfer operator.

\smallskip

For $g\in C(X)$, \ $g\geq 0$, we denote by $\ell(g)$ the logarithm of the spectral radius of $Ag$.

\begin{remark} \label{r4a-1}
If $g > 0$, then $Ag = A_{\ln g}$. In addition,
\begin{equation} \label{e-l-ell}
 \ell(g) = \lambda(\ln g), \qquad g>0.
\end{equation}
\end{remark}

Recall once more that whenever we use a transfer operator $A$ and the essential set $X_\alpha$
we assume that $A$ is $X_\alpha$-compatible and $\alpha$ is a finite-sheeted cover on $X_\alpha$.

\smallskip

The extension of Theorem~\ref{t-3-1a} on the situation in question is

\begin{theorem} \label{t-4-1a}
Let the cocycle\/ $\rho$ be continuous on\/ $X_\alpha$ $($in particular, this is true when\/ $\alpha \!: X_\alpha\to X_\alpha$ is a local homeomorphism\/$)$. Then for each\/ $g\in C(X)$, \ $g\geq 0$,
$$
 \ell(g) = \ell(\alpha,\rho g),
$$
$($here we recall that\/ $\ell(\alpha,a)$ in\/ \eqref{e.sp-exp} exploits only the values of\/ $a$ on\/ $X_\alpha$$)$.
\end{theorem}

\proof. Choose a sequence of strictly positive continuous functions $g_n \searrow\pin g$. By upper
semicontinuity of the spectral radius we have
\begin{equation} \label{et-4-1-0}
 \ell(g_n) \searrow\pin \ell(g).
\end{equation}
By Theorem~\ref{t-3-1a} and \eqref{e-l-ell} one obtains
\begin{equation} \label{et-4-1-1}
 \ell(g_n) = \ell(\alpha,\rho g_n).
\end{equation}
And from the explicit form of $\ell(\alpha,a)$ in \eqref{e.sp-exp} we conclude that
\begin{equation} \label{et-4-1-2}
 \ell(\alpha,\rho g_n) \searrow\pin \ell(\alpha, \rho g).
\end{equation}
Now \eqref{et-4-1-0}, \eqref{et-4-1-1}, and \eqref{et-4-1-2} imply
$$
\ell(g) = \ell(\alpha, \rho g). \qed
$$

The extension of Theorem~\ref{5..6} on the situation in question is

\begin{theorem} \label{t4a-1}
 {\rm (variational principle for transfer operators with nonnegative weights, see
 \cite[Theorem 11.2.]{AnBakhLeb11})}
Let\/ $Ag$ be a transfer operator defined in equation\/ \eqref{e4a-1}, where\/ $g\in C(X)$ and\/
$g\geq 0$. Then the following variational principle holds true\/$:$
\begin{equation} \label{e4a-2}
 \ell(g) \pin= \max_{\mu\in M_\alpha(X)}\biggl(\int_X \ln g\, d\mu +\tau(\mu) \biggr).
\end{equation}
\end{theorem}

Recalling Ruelle--Walters variational principle for topological pressure \eqref{e.VP-TP} one can
set for $a\in C(X)$, \ $a\geq 0$,
\begin{equation} \label{e4a-TP}
 P(\alpha,\ln a) := \sup_{\mu\in M_\alpha(X)} \biggl(\int_X \ln a\, d\mu +h_\alpha(\mu) \biggr).
\end{equation}

Now the extension of Theorem~\ref{t3-2} on the situation in question is

\begin{theorem} \label{t4a-2}
Let\/ $X$ be a compact metric space, $\alpha \!: X\to X$ be a non-contracting local homeomprphism on\/ $X_\alpha$, and\/ $g\in C(X)$, \ $g\geq 0$. Then for any nonnegative continuous extension of the cocycle\/ $\rho$ from\/ $X_\alpha$ onto\/ $X$ we have
$$
 \ell(g) = P\bigl(\alpha,\ln (g\rho)\bigr),
$$
where the topological pressure is defined by\/ \eqref{e4a-TP} $($we recall that\/ $\rho$ is continuous on\/ $X_\alpha$ by Corollary \ref{c-coc3}\pin$)$.
\end{theorem}

\proof. Let a transfer operator $A_{X_\alpha}\!: C(X_\alpha) \to C(X_\alpha)$ be given by
\eqref{e-transf-Y} with $Y=X_\alpha$ and denote by $\ell_{X_\alpha}(g)$ the logarithm of the
spectral radius of $A_{X_\alpha}g$. Recalling Theorems~\ref{t-tau-omega} and \ref{t-lambda-omega}
and exploiting Theorem~\ref{t4a-1} we conclude that
\begin{equation} \label{e4a-Omega}
 \ell (g) \pin=\pin \ell_{X_\alpha}(g) \pin= \max_{\mu\in M_\alpha(X)}
 \biggl(\int_{X_\alpha} \ln g\, d\mu +\tau_{X_\alpha}(\mu) \biggr).
\end{equation}
Choose strictly positive continuous functions $g_n \searrow\pin g$ and strictly positive continuous functions $\rho_n \searrow\pin \rho$. Let $A_n\!: C(X_\alpha) \to C(X_\alpha)$ be transfer operators associated with cocycles $\rho_n$ and consider the arising transfer operators $A_n g_n\!: C(X_\alpha) \to C(X_\alpha)$, and set $\ell_{X_\alpha}(g_n)$ to be the logarithm of the spectral radius of $A_n g_n$. By Theorem~\ref{t3-2} one has
$$
 \ell_{X_\alpha}(g_n) = P\bigl(\alpha, \ln(g_n\rho_n)\bigr).
$$
In addition, by upper semicontinuity of the spectral radius we obtain
\begin{equation} \label{e4a-5}
 \ell_{X_\alpha}(g_n) \searrow\pin \ell_{X_\alpha}(g),
\end{equation}
and we also have
\begin{equation} \label{e4a-6}
 P\bigl(\alpha, \ln (g_n\rho_n)\bigr)\searrow\pin P\bigl(\alpha, \ln(g\rho)\bigr),
\end{equation}

 \medskip\noindent
where the topological pressure is defined by \eqref{e4a-TP}. \qed

\medskip

For convenience of further reasoning we mention the next observation which naturally should be considered as a folklore.

\begin{lemma} \label{l4a-folk}
Let the entropy map\/ $M_\alpha(X) \ni \mu \mapsto h_\alpha(\mu) \in [0,\infty)$ be upper
semicontinuous\/ $($in $^*$-weak topology\/$)$ and\/ $g\in C(X)$, \ $g\geq 0$, then variational
principle \eqref{e4a-TP} reduces to
\begin{equation} \label{e4-h-upper}
 P(\alpha,\ln g) \pin= \max_{\mu\in M_\alpha(X)}
 \biggl(\int_{X_\alpha} \ln g\, d\mu +h_\alpha(\mu) \biggr).
\end{equation}
\end{lemma}

\proof. We have already mentioned that for $\mu \in M_\alpha(X)$ one has $\mathop{\mathrm{supp}}
\mu \subset X_\alpha$. Also the map $ M_\alpha(X) \ni \mu \mapsto \int_{X_\alpha} \ln g\, d\mu$ is
always upper semicontinuous. Therefore supremum in \eqref{e4a-TP} can be replaced by maximum in \eqref{e4-h-upper}. \qed

\medskip

Combining this lemma with Theorem~\ref{t4a-2} one gets

\begin{corollary} \label{c4-upper-semi}
Let\/ $X$ be a compact metric space and\/ $\alpha\! :X\to X$ be a non-contracting local homeomorphism of\/ $X_\alpha$ such that the entropy map\/ $M_\alpha(X) \ni \mu \mapsto
h_\alpha(\mu) \in [0,\infty)$ is upper semicontinuous. Then for each\/ $g\in C(X)$, \ $g\geq 0$ we
have
\begin{equation} \label{,,81}
 \ell(g) \pin= \max_{\mu\in M_\alpha(X)} \biggl(\int_{X_\alpha}\ln(g\rho)\,d\mu +h_\alpha(\mu)\biggr).
\end{equation}
\end{corollary}

\medskip

There is quite a number of dynamical systems on a metric space $(X, d)$ for which upper
semicontinuity of entropy map holds. Among them are, for example, expanding maps.
 Thus we also get the next

\begin{corollary} \label{c4-expand}
Let\/ $X$ be a compact metric space and $\alpha \!: X\to X$ be an expanding local homeomorphism on\/ $X_\alpha$. Then for each\/ $g\in C(X)$, \ $g\geq 0$, we have equality\/ \eqref{,,81}.
\end{corollary}

As is known $h_\alpha(\mu)$ is a concave function and therefore one can `close' this function to make it upper semicontinuous. Namely, set
\begin{equation} \label{e-h-reg}
 \bar{h}_\alpha(\mu):= \limsup_{\nu \to \mu} h_\alpha(\nu),
\end{equation}
where $\nu\to\mu$ is taken in $^*$-weak topology. This function $\bar{h}_\alpha(\mu)$ is concave
and upper semicontinuous on $M_\alpha(X)$. Replacing entropy $h_\alpha(\mu)$ by
$\bar{h}_\alpha(\mu)$ in the proof of Lemma~\ref{l4a-folk} one obtains

\begin{lemma} \label{l4a-folk1}
Let\/ $X$ be a compact metric space, $\alpha \!: X\to X$ be a continuous mapping and\/ $g\in C(X)$,
\ $g\geq 0$, then along with variational principle\/ \eqref{e4a-TP} we have
\begin{equation} \label{e4-h-upper1}
 P(\alpha,\ln g) \pin= \max_{\mu\in M_\alpha(X)}
 \biggl(\int_{X_\alpha} \ln g\, d\mu +\bar{h}_\alpha(\mu) \biggr).
\end{equation}
\end{lemma}

And as an analogue of Corollary~\ref{c4-upper-semi} for the `closed' entropy
$\overline{h}_\alpha(\mu)$ one has

\begin{corollary} \label{c4-upper-semi1}
Let\/ $X$ be a compact metric space and $\alpha\!: X\to X$ be a non-contracting local homeomorphism on\/ $X_\alpha$. Then for each\/ $g\in C(X)$, \ $g\geq 0$, we have
$$
 \ell(g) \pin= \max_{\mu\in M_\alpha(X)}
 \biggl(\int_{X_\alpha} \ln (g\rho)\, d\mu +\bar{h}_\alpha(\mu) \biggr).
$$
\end{corollary}

\medskip

We finish this section with relating $\ell(g)$ with integrals. As an immediate consequence of
Theorems~\ref{t.VP-lambda}, \ref{t2-invert}, \ref{t-4-1a}, and Corollary~\ref{c-coc3} one obtains

\begin{theorem} \label{t3-3}
Let the inverse rami-rate be zero\/ $(\omega(\alpha)=0),$ and the cocyle\/ $\rho$ on $X_\alpha$ be  continuous {\rm(}in particular, this takes place when\/ $\alpha$ is a homeomorphism on\/ $X_\alpha${\rm)}. Then
\begin{equation} \label{e3-erg}
 \ell(g) \pin= \max_{\mu\in E\!\pin M_\alpha(X)} \int_{X_\alpha} \ln(g\rho) \, d\mu \pin=
 \max_{\mu\in M_\alpha(X)} \int_{X_\alpha} \ln (g\rho) \, d\mu\,.
\end{equation}
\end{theorem}

\begin{remark} \label{r3-erg}
If $\alpha\!: X \to X$ is a homeomorphism then the description of transfer operators given in Subsection~\ref{ss-11} implies
$$
 [Af](x) = \rho\big(\alpha^{-1}(x)\big)f\big(\alpha^{-1}(x)\big),
$$
where $\rho \in C(X)$ is a nonnegative function. That is, $A$ is a weighted shift operator and so
also is the operator
$$
 [(Ag) f] (x) = [\rho\pin gf]\big(\alpha^{-1}(x)\big).
$$

 \medskip\noindent
Variational principles of \eqref{e3-erg} type for abstract weighted shift operators associated with
commutative Banach algebras automorphisms generated by isometries where worked out in
\cite{Kitover} and \cite{Lebedev79} (see also \cite[4]{Anton_Lebed} and \cite[5]{Anton}).
A comprehensive analysis of the corresponding variational principles and their interrelations as with integrals so also with Lyapunov exponents for abstract weighted shift operators associated with endomorphisms of Banach algebras is presented in \cite{KL20}.
\end{remark}

\section[$\boldsymbol T$-entropy vs Kolmogorov--Sinaj entropy and integrals]
         {$\boldsymbol T$-entropy vs Kolmogorov--Sinaj entropy and\\ integrals} \label{s-4}

In the previous sections we analyzed relationships between spectral potential, topological pressure
and integrals with respect to invariant measures. The results obtained naturally give us an opportunity to analyze relationships between $t$-entropy, entropy and integrals. This is the theme of the present section.

Henceforth we assume that $A$ is a given transfer operator for a dynamical system $(X,\alpha)$,
provided $A$ is $X_\alpha$-compatible and $\alpha$ is a finite-sheeted cover on $X_\alpha$; and $\rho$ is a cocycle on $X_\alpha$ defined by~\eqref{e-3.1}.

We recall one more observation that will be exploited in sequel.

In \cite[Propositions~8.4, 8.6]{AnBakhLeb11} it is proven that $t$-entropy map $\mu\mapsto
\tau(\mu)$ is a concave and upper semicontinuous function (in $^*$-weak topology). Therefore
formula \eqref{5,,12}:
$$
 \lambda(\psi) \pin=\max_{\mu\in M_\alpha(X)} \bigl(\mu[\psi] + \tau(\mu) \bigr)
$$
means that the spectral potential $\lambda(\psi)$ is nothing else than the Fenchel--Legendre
transform of $-\tau(\mu)$. Moreover, by the Fenchel--Legendre--Moreau duality upper
semicontin\-ui\-ty of $\tau(\mu)$ implies the equality
\begin{equation} \label{,,84}
 -\tau(\mu) \,=\, \inf_{\psi\in C(X,\mathbb R)} \big(\mu[\psi] -\lambda(\psi)\big),
\end{equation}
which means that $-\tau(\mu)$ is the Fenchel--Legendre dual functional to $\lambda(\psi)$. Therefore $\tau(\mu)$ is uniquely defined by the spectral potential $\lambda(\psi)$. By the
mentioned Fenchel--Legendre--Moreau duality we also conclude that if $S(\mu)$ is a certain concave
and upper semicontinuous (in $^*$-weak topology) function of $\mu$ such that
\begin{equation}  \label{e5-1}
 \lambda(\psi) \pin= \sup_{\mu\in M_\alpha(X)} \bigl(\mu[\psi]+S(\mu)\bigr)
\end{equation}
(i.\,e., $\lambda(\psi)$ is the Fenchel--Legendre transform of $-S(\mu)$) then
\begin{equation}  \label{e5-2}
 S(\mu) = \tau(\mu), \qquad \mu\in M_\alpha(X).
\end{equation}

Recall once more that whenever we use a transfer operator $A$ and the essential set $X_\alpha$
we assume that $A$ is $X_\alpha$-compatible.

Our first observation is relationship between $t$-entropy and integrals.

\begin{theorem} \label{t4-invert}
Let the inverse rami-rate be zero\/ $(\omega(\alpha)=0)$ and the cocycle\/ $\rho$ on $X_\alpha$ be  continuous\/ {\rm(}in particular, this takes place when\/ $\alpha$ is a homeomorphism on\/
$X_\alpha${\rm)}. Then for\/ $\mu\in M_\alpha(X)$ we have
\begin{equation} \label{e4-erg}
 \tau(\mu) \pin= \int_{X_\alpha}\ln\rho\, d\mu\pin.
\end{equation}
\end{theorem}

 \medskip

\proof. By Theorem~\ref{t3-3} for $\psi\in C(X)$ one has
\begin{equation} \label{e4-3}
 \lambda (\psi)\pin =\max_{\mu\in M_\alpha(X)}\biggl(\mu[\psi] +\int_{X_\alpha}\ln\rho\,d\mu\biggr).
\end{equation}
Note that the function
$$
M_\alpha(X) \ni \mu \longmapsto \int_{X_\alpha}\ln \rho\, d\mu
$$

 \medskip\noindent
is linear and upper semicontinuous. This observation along with \eqref{e5-1} and \eqref{e5-2}
implies the assertion of the theorem. \qed

\begin{remark} \label{r4-hom}
In the case when $\alpha\!: X\to X$ is a homeomorphism and $Af(x) =f\bigl(\alpha^{-1}(x)\bigr)$
the corresponding formula for $t$-entropy was obtained in \cite{BKKL19}.
\end{remark}

The next observation links $t$-entropy with Kolmogorov--Sinaj entropy.

\begin{theorem} \label{t4-t-entr-entrop}
Let\/ $X$ be a compact metric space, $\alpha$ be an open and non-contracting on\/~$X_\alpha$, and the entropy map\/ $M_\alpha(X) \ni \mu \mapsto h_\alpha(\mu) \in [0,\infty)$ be upper semicontinuous\/ $($in particular, this takes place when\/ $\alpha$ is open and expanding on\/ $X_\alpha$$)$. Then for\/ $\mu\in M_\alpha(X)$ we have
\begin{equation} \label{e4-entr}
 \tau(\mu) \pin= \int_{X_\alpha}\ln \rho\, d\mu + h_\alpha(\mu).
\end{equation}
\end{theorem}

\medskip

\proof. By Corollaries~\ref{c4-upper-semi}, \ref{c4-expand}, and \eqref{e-l-ell}, for $\psi\in C(X)$ one has
$$
 \lambda(\psi) \pin= \max_{\mu\in M_\alpha(X)} \biggl(\mu[\psi] +
 \biggl[\int_{X_\alpha}\ln\rho\,d\mu + h_\alpha(\mu)\biggr]\biggr).
$$
And by the condition of the theorem the function
$$
 M_\alpha(X) \ni \mu \longmapsto \int_{X_\alpha}\ln \rho\, d\mu + h_\alpha(\mu)
$$
is upper semicontinuous. Now the assertion of the theorem follows from \eqref{e5-1} and
\eqref{e5-2}. \qed

\begin{remark} \label{r4-1}
In the situation when $\alpha\!: X\to X$ is open and expanding formula \eqref{e4-entr}
for $t$-entropy was obtained in \cite{BKKL19} and \cite{BK19}.
\end{remark}

Replacing in the proof of the previous theorem entropy $h_\alpha(\mu)$ by its `closure'
$\bar{h}_\alpha(\mu)$ \eqref{e-h-reg} and exploiting Corollary~\ref{c4-upper-semi1} one obtains

\begin{theorem} \label{t4-t-entr-entrop1}
Let\/ $X$ be a compact metric space, and\/ $\alpha$ be open and non-contracting on\/ $X_\alpha$. Then for\/ $\mu\in M_\alpha(X)$ we have
$$
 \tau (\mu) \pin= \int_{X_\alpha}\ln\rho\, d\mu + \bar{h}_\alpha(\mu).
$$
\end{theorem}

\begin{remark} \label{r4-2}
In all the statements the set $X_\alpha$ and the assumption that $A$ is $X_\alpha$-compatible are essential. In fact the properties of $X_\alpha$ that were exploited are the following:

\smallskip

1) each $\mu \in M_\alpha (X)$ is supported on $X_\alpha$,

\smallskip

2) the operator $A$ is compatible with this set.

 \smallskip \noindent
Any set $Y\subset X$ possessing these two properties can be exploited in all the statements. For
example, one can take the set $\Omega(\alpha)$ of non-wandering points when $A$  is $\Omega(\alpha)$-compatible.
Or simply take the whole $X$.

In fact, the essential set $X_\alpha$ is the minimal one possessing the mentioned properties
whenever $A$ is compatible with it.
\end{remark}


\end{document}